\documentclass[12pt]{amsart}

\newtheorem{lem}{Lemma}[section]
\newtheorem{defi}[lem]{Definition}
\newtheorem{prop}[lem]{Proposition}
\newtheorem{rema}[lem]{Remark}
\newtheorem{theo}[lem]{Theorem}
\newtheorem{exam}[lem]{Example}

\usepackage{epsfig}

\begin{document}
\parindent0em
\title{Gradient Flows of closed 1-forms and their closed orbits}
\author{D. Sch\"utz}
\address{Department of Mathematics, SUNY Binghamton, Binghamton, NY 13902-6000}
\email{dirk@math.binghamton.edu}
\subjclass{Primary 57R70; Secondary 57R25}
\keywords{Novikov complex, closed 1-forms, zeta function, Dennis trace}
\begin{abstract}
In \cite{pajirn,pajitn}, Pajitnov considers the closed orbit structure of generic gradient flows of
circle-valued Morse functions. It turns out that the torsion of a chain homotopy equivalence
between the Novikov complex and the completed simplicial chain complex of the universal cover detects
the eta function of the flow. This eta function counts the closed orbits and reduces to the
logarithm of the zeta function after abelianizing. We extend this result to the case of
closed 1-forms which are Morse. To relate the torsion to the eta function we use the Dennis trace.
\end{abstract}
\maketitle
\section{Introduction}
Given a vector field on a smooth closed manifold $M$ there is a corresponding dynamical system and one
can investigate the closed orbits of this flow. It is desirable to collect all closed orbits in one power
series and study the algebraic topology and $K$-theory of this object. To do this observe that closed orbits represent
elements in $H_1(M)$ and also in the set of conjugacy classes of $\pi_1(M)$. We set $G=\pi_1(M)$.\\[0.2cm]
In \cite{fried}, Fried defines a commutative zeta function for certain nonsingular flows as a formal
power series and relates it to a Reidemeister torsion invariant of the manifold.\\[0.2cm]
The first noncommutative invariant for flows was introduced in Geoghegan and Nicas \cite{geonic} for
suspension flows. Their analogue of a zeta function is what they call the Lefschetz-Nielsen series which
lives in an infinite product of $0$-dimensional Hochschild homology groups.\\[0.2cm]
In the case of vector fields with singularities the first papers to obtain relations between zeta functions
and torsion are Hutchings and Lee \cite{hutlee,hutle2} and Pajitnov \cite{pajirn}, both dealing with gradients
of circle-valued Morse functions and with commutative invariants. In that situation the torsion invariant
no longer depends only on the topology of $M$ but the critical points enter via the Novikov complex.
Both papers have been generalized, Hutchings \cite{hutcth,hutchi} discusses closed 1-forms, still in
a commutative setting, while Pajitnov \cite{pajitn} gets a noncommutative result for circle-valued Morse
functions.\\[0.2cm]
Circle-valued Morse functions correspond to closed rational Morse 1-forms.
This paper discusses the noncommutative case for arbitrary closed Morse 1-forms. The geometric methods largely follow
Pajitnov \cite{pajitn}. In fact, the geometry in \cite{pajitn} is mainly contained in his earlier paper
\cite{pajirn}. The main difficulty is that the algebra required to keep track of the invariants is
more complicated than in the commutative case.
So instead of looking at a zeta function, Pajitnov \cite{pajitn} and we look at an eta function (or pre-zeta
function, compare Fried \cite[\S 2]{fried}), which generalizes the logarithm of the zeta function of the commutative
case. Since the conjugacy classes of $G$ do not form a group, we cannot take the exponential
function of this eta function. To compare this eta function with a certain torsion one needs a
logarithm-like homomorphism $\mathfrak{L}$ from $K_1$ of the Novikov ring to the object containing the eta function.
We depart somewhat in the definition of $\mathfrak{L}$ from Pajitnov \cite{pajitn} in that we take
a detour through Hochschild homology using the Dennis trace, compare Geoghegan and Nicas \cite[\S 5]{geonic}.
The main theorem we get is
\begin{theo}Let $\omega$ be a closed Morse 1-form on a smooth connected closed manifold $M^n$. Let $\xi:G\to\mathbb{R}$ be induced by $\omega$ and let $C^\Delta_\ast(\tilde{M})$ be the simplicial $\mathbb{Z}G$ complex coming from a smooth triangulation of $M$. For every $v\in\mathcal{G}_0(\omega)$ there is a natural chain homotopy equivalence $\varphi(v):\widehat{\mathbb{Z}G}_\xi\otimes_{\mathbb{Z}G}C_\ast^\Delta(\tilde{M})\to C_\ast(\omega,v)$ whose torsion $\tau(\varphi(v))$
lies in $\overline{W}$ and satisfies
\[\mathfrak{L}(\tau(\varphi(v)))=\eta(-v).\]
\label{intro}
\end{theo}
This theorem was obtained by Pajitnov in \cite{pajitn} in the rational case.
Here $v$ is the vector field whose eta function we look at, $C_\ast(\omega,v)$ is the
Novikov complex, a complex over the Novikov ring $\widehat{\mathbb{Z}G}_\xi$ and $\overline{W}$ a
particular subgroup of $\overline{K}_1^G(\widehat{\mathbb{Z}G}_\xi)$.
The set $\mathcal{G}_0(\omega)$ is a set of $C^0$-generic vector fields which are gradient with
respect to $\omega$, see Section 4.
The chain homotopy equivalence can be described as follows: given a smooth triangulation of $M$, we
can adjust this triangulation so that each simplex is transverse to the unstable manifolds of the
critical points of $\omega$. Then for a $k$-simplex $\sigma$ we define
\begin{equation}\label{formula}
\varphi(v)(\sigma)=\sum_{p\in{\rm crit}_k(\omega)}[\sigma:p]\,p
\end{equation}
where ${\rm crit}_k(\omega)$ is the set of critical points of $\omega$ having index $k$ and $[
\sigma:p]\in\widehat{\mathbb{Z}G}_\xi$ is the intersection number of a lifting of $\sigma$ to
$\tilde{M}$ with translates of the unstable manifold of a lifting of the critical point $p$. This
chain homotopy equivalence is basically described in Hutchings and Lee \cite[\S 2.3]{hutlee}. The
restriction that $v$ lie in $\mathcal{G}_0(\omega)$, a geometric condition due to Pajitnov \cite{pajirn},
then allows us to identify the torsion of $\varphi(v)$. This is achieved using the work of Farber
and Ranicki \cite{farran} and Ranicki \cite{ranick}. We choose a triangulation of $M$ such that
$\varphi(v)$ factors through a complex $\widehat{\mathbb{Z}G}_\xi\otimes C(v)_\ast$, where
$C(v)_\ast$ is a $\mathbb{Z}G$ complex which comes from a handlebody decomposition
on a codimension 1 submanifold $N$ that separates $M$ and a handlebody decomposition on the cobordism obtained by splitting along $N$.
It turns out that the complex $C(v)_\ast$ is the mapping cone of an injective $\mathbb{Z}G$ homomorphism
which depends on the vector field $v$. After
tensoring with the Novikov ring the natural projection to the cokernel is a chain homotopy equivalence.
But for $v\in\mathcal{G}_0(\omega)$ the cokernel can be identified with the Novikov complex.\\[0.2cm]
If the Novikov complex is not acyclic the torsion of a chain homotopy equivalence is not determined
by the complexes and we will give an example of two $\omega$-gradients
$v,w$ with $C_\ast(\omega,v)=C_\ast(\omega,w)$, but $\tau(\varphi(v))\not=\tau(\varphi(w))$,
see Remark \ref{rm54}.\\[0.2cm]
As mentioned before this paper is closely related to Pajitnov \cite{pajirn,pajitn}. The work of Hutchings \cite{hutcth,hutchi}
and Hutchings and Lee \cite{hutlee,hutle2} is in the same spirit, but with quite different methods. In
particular, Hutchings \cite{hutchi} contains a Theorem (see Theorem \ref{htheo} for the precise
statement) which might be considered a commutative version of Theorem \ref{intro}. The role of $\tau(\varphi(v))$
is played by two Reidemeister torsions. We show in section 5 how to recover Hutchings' theorem for
vector fields in $\mathcal{G}_0(\omega)$ as a corollary of Theorem \ref{intro}. In fact we obtain a stronger
``commutative theorem''; see Example \ref{examp}.\\[0.2cm]
The author would like to thank the referee for several suggestions, in particular for formula
(\ref{formula}). This paper will form a portion of the author's doctoral dissertation
written at the State University of New York at Binghamton under the direction of Ross Geoghegan.
\section{Morse theory of closed 1-forms}
\subsection{Novikov Rings}
Let $G$ be a group and $\xi:G\to\mathbb{R}$ be a homomorphism. For a ring $R$ we denote by $\widehat{\widehat{RG}}$ the abelian group of all functions $G\to R$. For $\lambda\in\widehat{\widehat{RG}}$ let \\supp $\lambda=\{g\in G\,|\,\lambda(g)\not=0\}$. Then we define
\[\widehat{RG}_\xi=\{\lambda\in\widehat{\widehat{RG}}\,|\,\forall r\in\mathbb{R}\hspace{0.4cm}\#\,\mbox{supp }\lambda\cap\xi^{-1}([r,\infty))<\infty\}\]
For $\lambda_1,\lambda_2\in\widehat{RG}_\xi$ we set $(\lambda_1\cdot\lambda_2)(g)=\sum\limits_{h_1,h_2\in G\atop h_1h_2=g}\lambda_1(h_1)\lambda_2(h_2)$, then $\lambda_1\cdot\lambda_2$ is a well defined element of $\widehat{RG}_\xi$ and turns $\widehat{RG}_\xi$ into a ring, the \em Novikov ring\em. It contains the usual group ring $RG$ as a subring and we have $RG=\widehat{RG}_\xi$ if and only if $\xi$ is the zero homomorphism.
\begin{defi}\em The \em norm \em of $\lambda\in\widehat{RG}_\xi$ is defined to be
\[\|\lambda\|=\|\lambda\|_\xi=\inf\{t\in(0,\infty)|\mbox{ supp }\lambda\subset\xi^{-1}((-\infty,\log t])\}\]
\end{defi}
It has the following nice properties:
\begin{enumerate}
\item $\|\lambda\|\geq 0$ and $\|\lambda\|=0$ if and only if $\lambda=0$.
\item $\|\lambda\|=\|-\lambda\|$.
\item $\|\lambda+\mu\|\leq\max\{\|\lambda\|,\|\mu\|\}$.
\item $\|\lambda\cdot\mu\|\leq\|\lambda\|\cdot\|\mu\|$.
\end{enumerate}
If $N$ is a normal subgroup of $G$ that is contained in $\ker\xi$ we get a well defined homomorphism $\bar{\xi}:G/N\to\mathbb{R}$ and a well defined ring epimorphism $\varepsilon:\widehat{RG}_\xi\to\widehat{RG/N}_{\bar{\xi}}$ given by \\$\varepsilon(\lambda)(gN)=\sum\limits_{n\in N}\lambda(gn)$.\\[0.2cm]
Now let $\Gamma$ be the set of conjugacy classes of $G$. Again the homomorphism $\xi$ induces a well defined map $\Gamma\to\mathbb{R}$ which we also denote by $\xi$. In analogy with above we define $\widehat{R\Gamma}_\xi$, but since there is no well defined multiplication in $\Gamma$, this object is just an abelian group. Again there is an epimorphism $\varepsilon:\widehat{RG}_\xi\to\widehat{R\Gamma}_\xi$ of abelian groups. We can think of $\widehat{R\Gamma}_\xi$ as lying between $\widehat{RG}_\xi$ and $\widehat{RH_1(G)}_{\bar{\xi}}$. If $g\in G$, we denote the conjugacy class of $g$ by $\{g\}$.\\[0.2cm]
Now we will turn our attention to $K_1(\widehat{\mathbb{Z}G}_\xi)$. For the definition of $K_1$ we refer the reader to Cohen \cite{cohen} or Milnor \cite{milntr}.
First we disregard units of the form $\pm g$, hence look at $\overline{K}_1^G(\widehat{\mathbb{Z}G}_\xi)=K_1(\widehat{\mathbb{Z}G}_\xi)/\langle[\pm g]\rangle$. There is another type of ``elementary unit'' in $\widehat{\mathbb{Z}G}_\xi$, namely, let $a\in\widehat{\mathbb{Z}G}_\xi$ satisfy $\|a\|<1$. Then $\sum_{n=0}^\infty a^n$ is a well defined element of $\widehat{\mathbb{Z}G}_\xi$ and the inverse of $1-a$.
These form a subgroup of the units in  $\widehat{\mathbb{Z}G}_\xi$. We denote the image of this subgroup in $\overline{K}_1^G(\widehat{\mathbb{Z}G}_\xi)$ by $\overline{W}$.
\subsection{Closed 1-forms and Vector Fields}
Let $M^n$ be a closed connected smooth manifold. By de Rham's theorem $\{$closed $1$-forms on $M\}/\{$exact $1$-forms on $M\}\cong H^1(M;\mathbb{R})\cong\\ $Hom$(H_1(M),\mathbb{R})$, so a closed $1$-form $\omega$ induces a homomorphism $\xi_\omega:\pi_1(M)\to\mathbb{R}$ which can be explicitely stated by the formula $\xi_\omega(g)=\int_\gamma\omega\in\mathbb{R}$, where $\gamma$ is a smooth loop representing $g\in\pi_1(M)$. 
Set $G=\pi_1(M)$. Then $G$ is finitely presented, so the image of $\xi_\omega$ is a finitely generated subgroup of $\mathbb{R}$, hence isomorphic to $\mathbb{Z}^k$ for some integer $k$. If $k=1$ $\omega$ is said to be \em rational\em, if $k>1$ it is \em irrational\em.\\[0.2cm]
Rational $1$-forms can be described by circle valued functions $f:M\to S^1$ in the following way: Let $p:\mathbb{R}\to\mathbb{R}/\mathbb{Z}=S^1$ be the usual covering projection, let $\alpha$ be the closed $1$-form on $S^1$ such that $p^*\alpha=dx$; then $f^*\alpha$ is a closed $1$-form and im $\xi_{f^*\alpha}\subset\mathbb{Z}\subset\mathbb{R}$. To obtain other infinite cyclic subgroups of $\mathbb{R}$ as images of $\xi$ one uses circles of different size.\\[0.2cm]
Now, given a rational $1$-form $\omega$ there is an infinite cyclic covering space $q:\bar{M}\to M$ such that $q^*\omega=d\bar{f}$, namely the one corresponding to $\ker\xi_\omega$. Let $t$ be the generator of the covering transformation group of $\bar{M}$ with $\bar{f}(tx)>\bar{f}(x)$ for $x\in\bar{M}$. Then $\bar{f}$ defines a map $f:M\to\mathbb{R}/(\bar{f}(tx)-\bar{f}(x))\mathbb{Z}=S^1$ which induces a surjection on fundamental group.\\[0.2cm]
Notice that for irrational closed $1$-forms $\omega$ there is a $\mathbb{Z}^k$-covering space $q:\bar{M}\to M$ such that $q^*\omega=d\bar{f}$.\\[0.2cm]
Locally a closed $1$-form is exact. We will call a closed $1$-form a \em Morse form \em if $\omega$ is locally represented by the differential of real valued functions whose critical points are nondegenerate.
So if $\omega$ is a Morse form, then $\omega$ has only finitely many critical points and every critical point has a well defined index.
\begin{defi}\em Let $\omega$ be a closed $1$-form. A vector field $v$ is called an \em $\omega$-gradient\em, if there exists a Riemannian metric $g$ such that $\omega_x(X)=g(X,v(x))$ for every $x\in M$ and $X\in T_x M$.
\end{defi}
The next Lemma allows us to forget about the Riemannian metric and will be useful in using vector fields as gradients of different Morse forms.
\begin{lem}\label{smalo}Let $\omega$ be a Morse form and $v$ a vector field. Then $v$ is an $\omega$-gradient if and only if
\begin{enumerate}
\item For every critical point $p$ of $\omega$ there exists a neighborhood $U_p$ of $p$ and a Riemannian metric $g$ on $U_p$ such that $\omega_x(X)=g(X,v(x))$ for every $x\in U_p$ and $X\in T_x U_p$.
\item If $\omega_x\not=0$, then $\omega_x(v(x))>0$.
\end{enumerate}
\end{lem}
\begin{proof}
The ``only if'' direction is clear.
For the ``if'' direction choose disjoint neighborhoods $U_1,\ldots ,U_k$, each with a Riemannian metric coming from 1. for every critical point of $\omega$. Now choose finitely many contractible open sets $V_1,\ldots\hspace{-1pt},V_m$ with $\bigcup V_i\subset M-\{$critical points$\}$ that together with the $U_j$'s cover $M$. Using 2., it is easy to find a Riemannian metric on each $V_i$ that turns $v|_{V_i}$ into a gradient of $\omega|_{V_i}$. Now the required Riemannian metric is obtained by using a partition of unity.
\end{proof}
\begin{rema}\em
Some authors (e.g.\ Milnor \cite{milnhc}, Pajitnov \cite{pajito,pajirn}) use a more restricted
version for an $\omega$-gradient, namely, a sharper version of 1.\ in the Lemma. For an even more
general definition of $\omega$-gradient we refer the reader to Pajitnov \cite{pajiov}, which
contains most of the modifications on a vector field that we will need.\em
\end{rema}
\subsection{The Novikov Complex of a Morse form}
Given a Morse form $\omega$ and an $\omega$-gradient $v$ we denote for a critical point $p$ of $\omega$ the unstable, resp. stable, manifold of $p$ by $D_R(p)$, resp. $D_L(p)$. So if $\Phi:M\times\mathbb{R}\to M$ denotes the flow of $v$, then $D_R(p)=\{x\in M|\Phi(x,t)\to p$ for $t\to-\infty\}$ and $D_L(p)=\{x\in M|\Phi(x,t)\to p$ for $t\to\infty\}$.
If the index of $p$ is $i$, then $D_R(p)$ is an immersed open disk of dimension $n-i$ and $D_L(p)$ of dimension $i$. We say $v$ satisfies the \em transversality condition \em if all discs $D_L(p)$ and $D_R(q)$ intersect transversely for all critical points $p,q$ of $\omega$.\\[0.2cm]
Given a Morse form $\omega$ and an $\omega$-gradient $v$ satisfying the transversality condition we can define the \em Novikov complex \em $C_\ast(\omega,v)$ which is in each dimension $i$ a free $\widehat{\mathbb{Z}G}_\xi$ complex with one generator for every critical point of index $i$. Here $\xi$ is the homomorphism induced by $\omega$. The boundary homomorphism of $C_\ast(\omega,v)$ is based on the number of trajectories between critical points of adjacent indices. For more details see Pajitnov \cite{pajito} or Latour \cite{latour}.
This chain complex is chain homotopy equivalent to $\widehat{\mathbb{Z}G}_\xi\otimes_{\mathbb{Z}G}C_\ast^\Delta(\tilde{M})$, where $C_\ast^\Delta(\tilde{M})$ is the simplicial chain complex of the universal cover $\tilde{M}$ of $M$ with respect to a smooth triangulation of $M$ lifted to $\tilde{M}$.\footnote{Of course there are other Novikov complexes corresponding to other regular coverings of $M$ but we are mainly interested in the universal covering.}
Furthermore there is a chain homotopy equivalence whose torsion is in $\overline{W}\subset\overline{K}_1^G(\widehat{\mathbb{Z}G}_\xi)$. In the rational case this is proven in Pajitnov \cite{pajito}, for the general case see Latour \cite{latour}.
The map described in (\ref{formula}) in the introduction can be used for this. We will show this in
\ref{sb4.1} at least for an $\omega$-gradient $v$ satisfying a ``cellularity condition''.\\[0.2cm]
Let us discuss this map. A smooth triangulation $\Delta$ of $M$ is called \em adjusted to \em $v$, if
every $k$-simplex $\sigma$ intersects the unstable manifolds $D_R(p)$ transversely for all critical
points $p$ of index $\geq k$. To see the existence, assume $\omega$ is rational, the general case
follows by approximation, see \ref{sb4.2}.
A triangulation $\Delta$ lifts to a triangulation of $\bar{M}$, an
infinite cyclic covering space, compare \ref{CHtyTp}. For a diffeomorphism $\psi$ of $M$ we denote
by $\psi\Delta$ the triangulation of $M$ where simplices are composed with $\psi$. If we change the
triangulation of $M$ by an isotopy, we can get transverse intersections in $\bar{M}$ of lifted
simplices with finitely many unstable manifolds by the results of \ref{app1} in the appendix. Since
the results there give openness and density among diffeomorphisms we get a generic set of
diffeomorphisms $\psi$ of $M$ isotopic to the identity such that $\psi\Delta$ is adjusted
to $v$.\\[0.2cm]
Given an adjusted triangulation $\Delta$ we get a chain map
\[\varphi(v):\widehat{\mathbb{Z}G}_\xi\otimes_{\mathbb{Z}G}C^\Delta_\ast(\tilde{M})\to
C_\ast(\omega,v)\]
by formula (\ref{formula}). That $\varphi(v)$ is indeed a chain map follows from the exact case,
which is described in the appendix. Lemma \ref{alem2} also carries over so that different adjustments
to a triangulation lead to chain homotopic maps. Finally, if $\Delta'$ is a subdivision of $\Delta$
such that $\psi\Delta'$ is adjusted to $v$, so is $\psi\Delta$ and the diagram
\[
\begin{array}{rcl}\widehat{\mathbb{Z}G}_\xi\otimes_{\mathbb{Z}G}C_\ast^{\psi\Delta}
(\tilde{W},\tilde{M}_0)&\stackrel{\rm sd}{\longrightarrow}&
\widehat{\mathbb{Z}G}_\xi\otimes_{\mathbb{Z}G}C_\ast^{\psi\Delta'}(\tilde{W},\tilde{M}_0)\\[0.2cm]
\varphi(v)\searrow& &\swarrow\varphi(v)\\[0.2cm]
\multicolumn{3}{c}{C_\ast(\omega,v)}\end{array}\]
commutes. So once we show that $\varphi(v)$ is a chain homotopy equivalence, its torsion
does not depend on the triangulation.
\subsection{The Chain Homotopy type of \boldmath$\cdot\otimes_{\mathbb{Z}G}C_\ast^\Delta(\tilde{M})$.}
\label{CHtyTp}
The following is a construction of Farber and Ranicki \cite{farran} written to fit our purposes,
compare also Pajitnov \cite[\S 7]{pajirn}.\\[0.2cm]
Given a circle valued Morse map $f:M\to S^1$ which induces a surjection on fundamental group we get
a lifting $\bar{f}:\bar{M}\to\mathbb{R}$ where $\bar{M}$ is an infinite cyclic covering space.
Assuming that $0\in\mathbb{R}$ is a regular value, set $N=\bar{f}^{-1}(0)$, $M_N=\bar{f}^{-1}([0,1])$.
We get a handle decomposition of the cobordism $(M_N;N,tN=\bar{f}^{-1}(1))$ from the Morse function
$\bar{f}|_{M_N}:M_N\to[0,1]$. By choosing a cell decomposition of $N$, Farber and Ranicki
\cite{farran} construct a finitely generated free $\mathbb{Z}G$ complex $C(v)_\ast$ homotopy equivalent to
$C_\ast^\Delta(\tilde{M})$. Let us recall the construction from \cite[\S 3]{farran}. Let $p:
\tilde{M}\to\bar{M}$ be the universal covering projection and let $\tilde{f}:\tilde{M}\to\mathbb{R}$
be $\bar{f}\circ p$. The cell decomposition of $N$ leads to a cell decomposition of $\tilde{f}^{-1}
(\mathbb{Z})$ and the resulting cell complex, denoted by $D$, is a finitely generated free $\mathbb{Z}G$
complex. Let $c_i(N)$ be the number of $i$-cells in $N$. Now $E$ is a finitely generated free
$\mathbb{Z}G$ complex containing $D$ as a subcomplex and the remaining generators correspond to
critical points of $\bar{f}|_{M_N}$.\\[0.2cm]
The inclusions $N\hookrightarrow M_N$ and $tN\hookrightarrow M_N$ induce $\mathbb{Z}G$ chain maps
$i:D\to E$ and $k:D\to E$ and since $D$ is a subcomplex, $i:D\to E$ is a split injection. Then
$C(v)_\ast=\mathcal{C}(i-k:D\to E)$ is the mapping cone of $i-k$.
In particular rank $C(v)_i=c_i(N)+ c_{i-1}(N)+\#$ critical points of $f$ having index $i$. In
\ref{sb4.1} we will use the geometry of Pajitnov \cite{pajirn} to get a more detailed version of
this chain complex.\\[0.2cm]
Let $R$ be a ring and $\eta:\mathbb{Z}G\to R$ a ring homomorphism such that \\
id$_R\otimes_{\mathbb{Z}G}$ proj$_D(i-k):R\otimes_{\mathbb{Z}G}D\to R\otimes_{\mathbb{Z}G}D$ is an
automorphism. Then by the Deformation Lemma of Farber and Ranicki \cite{farran}, see also Ranicki
\cite[Prop. 1.9]{ranick}, the chain complex \\$R\otimes_{\mathbb{Z}G}C(v)_\ast$ is chain homotopy equivalent to coker$($id$_R\otimes_{\mathbb{Z}G}(i-k))=:\hat{C}$, a finitely generated free $R$ complex with rank $\hat{C}_i=\#$ critical points of $f$ having index $i$.
In fact the chain equivalence is identified in \cite{ranick} to be the natural projection \\$p:R\otimes_{\mathbb{Z}G}\mathcal{C}(i-k)\to$ coker$($id$_R\otimes(i-k))$.\\[0.2cm]
So to use the Deformation Lemma one has to turn a certain square matrix $I-A$ representing proj$_D(i-k)$ over $\mathbb{Z}G$ into an invertible matrix over a ring $R$. The matrix $A$ can be chosen to satisfy $\|A_{ij}\|_\xi<1$ for every entry of $A$, where $\xi$ is induced by $f$. Obvious candidates for $R$ are the noncommutative localization used by Farber and Ranicki \cite{farran} and the Novikov ring $\widehat{\mathbb{Z}G}_\xi$. A not so obvious candidate is a Novikov ring $\widehat{\mathbb{Z}G}_{\xi'}$ where $\xi'$ is ``close'' to $\xi$; this will be discussed in section 4.
\begin{rema}\em
Farber \cite{farber} has extended the Deformation Lemma of \cite{farran} to the case of closed $1$-forms using a certain noncommutative localization.\em
\end{rema}
Furthermore, Ranicki \cite[Prop. 1.9]{ranick} contains the calculation of the torsion of the chain homotopy equivalence $p:R\otimes_{\mathbb{Z}G}C(v)_\ast\to\hat{C}$. It is given by
\begin{equation}\label{torsion}
\tau(p)=\sum_{i=0}^{n-1}(-1)^{i+1}\tau(I-A_i:R\otimes_{\mathbb{Z}G}D_i\to R\otimes_{\mathbb{Z}G}D_i)\in\overline{K}_1(R).\end{equation}
\subsection{The Eta Function of a Gradient}
Let $v$ be a vector field. By a closed orbit of $v$ we mean a nonconstant map $\gamma:S^1\to M$ with $\gamma'(x)=v(\gamma(x))$. The \em multiplicity \em $m(\gamma)$ is the largest positive integer $m$ such that $\gamma$ factors through an $m$-fold covering $S^1\to S^1$. 
We say two closed orbits are \em equivalent \em if they only differ by a rotation of $S^1$. We denote the set of equivalence classes by $Cl(v)$. 
Notice that $\gamma\in Cl(v)$ gives a well defined element $\{\gamma\}\in\Gamma$.
A closed orbit $\gamma$ is called \em nondegenerate \em if $\det(I-dP)\not=0$, where $P$ is a Poincar\'e map corresponding to $\gamma$.
In that case we define $\varepsilon(\gamma)\in\{1,-1\}$ to be the sign of $\det(I-dP)$.\footnote{i.e. $\varepsilon(\gamma)$ is the fixed point index of $P$ at the isolated fixed point coming from the closed orbit.}\\[0.2cm]
Now let $\omega$ be a Morse form. We denote by $\mathcal{G}(\omega)$ the set of all $\omega$-gradients that satisfy the transversality condition and whose closed orbits are nondegenerate.
For $v\in \mathcal{G}(\omega)$ we define the \em eta-function \em of $-v$ to be the element of $\widehat{\mathbb{Q}\Gamma}_\xi$ defined by 
\[\eta(-v)(\{g\})=\sum_{\gamma\in Cl(-v)\atop \{\gamma\}=\{g\}}\frac{\varepsilon(\gamma)}{m(\gamma)}\]
Again $\xi$ is induced by $\omega$. For the proof that $\eta(-v)$ is a well defined element of $\widehat{\mathbb{Q}\Gamma}_\xi$ we refer the reader to Hutchings \cite[\S 3.2]{hutcth}.
\section{Algebraic Constructions}
\subsection{Hochschild Homology} Let $R$ be a ring and $S$ an $R$-algebra. For an $S-S$ bimodule $M$ we define the \em Hochschild chain complex \em $(C_\ast(S,M),d)$ by $C_n(S,M)=S\otimes\ldots\otimes S\otimes M$, where the product contains $n$ copies of $S$ and the tensor products are taken over $R$. The boundary operator is given by
\begin{eqnarray*}
d(s_1\otimes\ldots\otimes s_n\otimes m)&=&s_2\otimes\ldots\otimes s_n\otimes ms_1\\
 & &+\sum_{i=1}^{n-1}(-1)^i s_1\otimes\ldots\otimes s_is_{i+1}\otimes\ldots\otimes s_n\otimes m\\
 & &+(-1)^n s_1\otimes\ldots\otimes s_{n-1}\otimes s_nm\end{eqnarray*}
The $n$-th Hochschild homology group of $S$ with coefficients in $M$ is denoted by $H\!H_n(S,M)$. If $M=S$ and the bimodule structure is given by ordinary multiplication we write $H\!H_\ast(S)$ instead of $H\!H_\ast(S,M)$. We will mainly be interested in the case where $R=\mathbb{Z}$, $S=M=\widehat{\mathbb{Z}G}_\xi$ is a Novikov ring and $n=1$. A useful observation is that $d(1\otimes 1\otimes x)=1\otimes x$ and hence classes represented by $1\otimes x$ are automatically $0$ in $H\!H_1(S,M)$.\\[0.2cm]
Given an $n\times k$ matrix $A=(A_{ij})$ over $S$ and an $k\times n$ matrix $B=(B_{ij})$ over $M$ we define an $n\times n$ matrix $A\otimes B$ with entries in $S\otimes M$ by setting $(A\otimes B)_{ij}=\sum_{l=1}^kA_{il}\otimes B_{lj}$. The trace of this matrix, trace $A\otimes B$, is given by $\sum_{l,m}A_{lm}\otimes B_{ml}$ and is an element of $C_1(S,M)$, it is a cycle if and only if trace$(AB)=$ trace$(BA)$. For more on Hochschild homology see Geoghegan and Nicas \cite{geonic} or Igusa \cite{igusa}.
\subsection{The homomorphism \boldmath$\mathfrak{L}$} For a ring $R$ with unit, which, in view of subsection 3.1, we can think of as a $\mathbb{Z}$-algebra, there is the Dennis trace homomorphism $DT:K_1(R)\to H\!H_1(R)$ defined as follows: If $\alpha\in K_1(R)$ is represented by the matrix $A$, then \\$DT(\alpha)=[$trace $(A\otimes A^{-1})]\in H\!H_1(R)$, see Igusa \cite[\S 1]{igusa}.
It is easy to see that the Dennis trace factors through $\overline{K}_1(R)=K_1(R)/\langle[-1]\rangle$.\\[0.2cm]
We want to define a homomorphism $\mathfrak{L}:\overline{W}\to\widehat{\mathbb{Q}\Gamma}_\xi$. To do this define\\ $m:C_1(\widehat{\mathbb{Z}G}_\xi,\widehat{\mathbb{Z}G}_\xi)=\widehat{\mathbb{Z}G}_\xi\otimes\widehat{\mathbb{Z}G}_\xi\to\widehat{\mathbb{R}\Gamma}_\xi$ by
\[ m(\lambda_1\otimes\lambda_2):\gamma\mapsto\left\{\begin{array}{cc}
\displaystyle\sum_{h_1,h_2\in G\atop \{h_1h_2\}=\gamma}\frac{\xi(h_1)}{\xi(\gamma)}\lambda_1(h_1)\lambda_2(h_2) & \mbox{if }\xi(\gamma)<0\\
0 & \mbox{if }\xi(\gamma)\geq 0\end{array}\right. \]
We can think of $m$ as a weighted combination of multiplication in $\widehat{\mathbb{Z}G}_\xi$ and the augmentation $\varepsilon:\widehat{\mathbb{Z}G}_\xi\to\widehat{\mathbb{Z}\Gamma}_\xi$.
\begin{lem} The homomorphism $m$ induces a homomorphism $\mu:H\!H_1(\widehat{\mathbb{Z}G}_\xi)\to\widehat{\mathbb{R}\Gamma}_\xi$.
\end{lem}
\begin{proof}
It is to be shown that $m$ vanishes on boundaries. So let $\gamma\in\Gamma$ satisfy $\xi(\gamma)<0$, then
\begin{eqnarray*}
& &m(\lambda_2\otimes\lambda_3\lambda_1-\lambda_1\lambda_2\otimes\lambda_3+\lambda_1\otimes\lambda_2\lambda_3)(\gamma)\\
& &=\hspace{0.3cm}\sum_{g,h\in G \atop \{gh\}=\gamma}\frac{\xi(g)}{\xi(\gamma)}\lambda_2(g)\lambda_3\lambda_1(h)-\sum_{g,h\in G \atop \{gh\}=\gamma}\frac{\xi(g)}{\xi(\gamma)}\lambda_1\lambda_2(g)\lambda_3(h)+\sum_{g,h\in G \atop \{gh\}=\gamma}\frac{\xi(g)}{\xi(\gamma)}\lambda_1(g)\lambda_2\lambda_3(h)\end{eqnarray*}
\begin{eqnarray*}
\hspace{0.3cm}&=&\frac{1}{\xi(\gamma)}\left(\sum_{g_2,g_3,g_1\in G \atop \{g_2g_3g_1\}=\gamma}\xi(g_2)\lambda_2(g_2)\lambda_3(g_3)\lambda_1(g_1)-\sum_{g_1,g_2,g_3\in G \atop \{g_1g_2g_3\}=\gamma}\xi(g_1g_2)\lambda_1(g_1)\lambda_2(g_2)\lambda_3(g_3)\right.\\
\hspace{0.3cm}& &\left.+\sum_{g_1,g_2,g_3\in G \atop \{g_1g_2g_3\}=\gamma}\xi(g_1)\lambda_1(g_1)\lambda_2(g_2)\lambda_3(g_3)\right)\\
&=&0
\end{eqnarray*}
since $\{g_2g_3g_1\}=\{g_1g_2g_3\}$ and $\xi$ is a homomorphism.
\end{proof}
For $g\in G$ we have $m(g\otimes g^{-1})=0$, therefore the composition $\mu\circ DT$ factors through $\overline{K}_1^G(\widehat{\mathbb{Z}G}_\xi)$, call this homomorphism $L:\overline{K}_1^G(\widehat{\mathbb{Z}G}_\xi)\to\widehat{\mathbb{R}\Gamma}_\xi$. We want to examine how this homomorphism behaves on $\overline{W}$. For future purposes it will be useful not just to look at $1\times 1$ matrices.
\begin{defi}\em An $n\times n$ matrix $A$ over $\widehat{\mathbb{Z}G}_\xi$ is called \em $\xi$-regular\em, if there exists $K<0$ such that $\|A_{i_1i_2}\|\cdot\|A_{i_2i_3}\|\cdots\|A_{i_mi_1}\|\leq\exp(Km)$ for all $m\geq 1$, $1\leq i_j\leq n$, $1\leq j\leq m$.
\end{defi}
For example, a matrix $A$ which satisfies $\|A_{ij}\| < 1$ for every entry is $\xi$-regular,
but $\xi$-regular matrices can have entries $A_{ij}$ with $\|A_{ij}\|\geq 1$. The motivation for
$\xi$-regular matrices comes from our approximation arguments, see Remark \ref{apprem}.
\begin{lem}
Let $A$ be a $\xi$-regular matrix. Then
\begin{enumerate}
\item The matrix $I-A$ is invertible over $\widehat{\mathbb{Z}G}_\xi$ and the inverse is given by $I+A+A^2+\ldots$.
\item Denote the image of $I-A$ in $\overline{K}_1^G(\widehat{\mathbb{Z}G}_\xi)$ by $\tau(I-A)$; then $\tau(I-A)\in\overline{W}$.
\end{enumerate}
\end{lem}
\begin{proof}1. We need to show that $I+A+A^2+\ldots$ is a well defined matrix over $\widehat{\mathbb{Z}G}_\xi$. Note that $(A^m)_{ij}=\displaystyle\sum_{i_1,\ldots\hspace{-1pt},i_{m-1}=1}^nA_{ii_1}A_{i_1i_2}\cdots A_{i_{m-1}j}$.
We will look at terms of the form \\$A_\ast=A_{i_1i_2}A_{i_2i_3}\cdots A_{i_{m-1}i_m}$ and get an estimate for $\|A_\ast\|$. The idea is to write $A_\ast$ as a word $C_1D_1\cdots C_lD_l$ where the length of the word $C_1\cdots C_l$ is smaller than $n$ and the words $D_j$ are of the form $A_{j_1j_2}\cdots A_{j_kj_1}$.\\[0.2cm]
So assume that $A_\ast=A_{i_1i_2}A_{i_2i_3}\cdots A_{i_{m-1}i_m}$ where $m>n+1$. Let $i_j$ be the first index whose value appears more than once. Since these numbers are between $1$ and $n$ we have $j\leq n$. Let $k$ be the largest number such that $i_k=i_j$, then
\[A_\ast=A_{i_1i_2}\cdots A_{i_{j-1}i_j}A_{i_ji_{j+1}}\cdots A_{i_{k-1}i_k}A_{i_ki_{k+1}}\cdots A_{i_{m-1}i_m}=B_1\cdots B_{j-1}D_1B_jA_{1\ast}\]
where $B_1=A_{i_1i_2},\ldots$, $B_{j-1}=A_{i_{j-1}i_j}$, $D_1=A_{i_ji_{j+1}}\cdots A_{i_{k-1}i_k}$, $B_j=A_{i_ki_{k+1}}$ and \\$A_{1\ast}=A_{i_{k+1}i_{k+2}}\cdots A_{i_{m-1}i_m}$.
Notice that $\|D_1\|\leq\exp(K(k-j))$. Now look at $A_{1\ast}$; among the indices $i_{k+1},\ldots$,$i_{m-1}$ are at most $n-j$ numbers; the numbers $i_1,\ldots$,$i_j$, which are all different, do not appear. If $m-1-k>n-j$, one of these numbers will appear more than once.
Let $i_{k+j_1}$ be the first such index and $i_{k+k_1}$ the last index equal to $i_{k+j_1}$. Again we get $j_1\leq n-j$, hence $j+j_1\leq n$. As above we can write
\[A_\ast=B_1\cdots B_{j-1}D_1B_j\cdots B_{j+j_1-1}D_2B_{j+j_1}A_{2\ast}\]
We continue this process until we get
\[A_\ast=B_1\cdots B_{j-1}D_1B_j\cdots B_{j+j_1-1}D_2B_{j+j_1}\cdots\cdots D_r\cdots B_l\]
with $l\leq n$. Since $A_\ast$ consists of $m-1$ letters we get $\|D_1\|\cdots \|D_r\|\leq\exp(K(m-1-n))$.
Let $M\in\mathbb{R}$ be a number such that $\|A_{st}\|\leq M$ for all $s,t$. Then $\|A_\ast\|\leq\exp(K(m-1-n))\cdot M^n$. Since $K(m-1-n)\to-\infty$ as $m\to\infty$, $\sum\limits_{m=0}^\infty A^m_{ij}$ is a well defined element of $\widehat{\mathbb{Z}G}_\xi$.\\[0.2cm]
2. The argument is the same as in Pajitnov \cite[Prop. 1.2]{pajitn}. Using elementary row reductions we obtain a matrix of the form
\[\left(\begin{array}{cccc}1-A_{11}&-A_{12}&\cdots&-A_{1n}\\
0\\
\vdots&\multicolumn{3}{c}{I-A'}\\
0\end{array}\right)\]
where $A'$ is an $n-1\times n-1$ matrix which is again $\xi$-regular with the same $K$. Induction gives the result.
\end{proof}
\begin{prop}\label{logar}Let $A$ be a $\xi$-regular matrix over $\widehat{\mathbb{Z}G}_\xi$, then
\[L(\tau(I-A))=-\varepsilon\left(\sum_{m=1}^\infty\frac{\mbox{\em trace \em}A^m}{m}\right)\]
In particular, $L$ induces a homomorphism $\mathfrak{L}:\overline{W}\to\widehat{\mathbb{Q}\Gamma}_\xi$.
\end{prop}
\begin{proof}
As before denote the image of $I-A$ in $\overline{K}_1(\widehat{\mathbb{Z}G}_\xi)$ by $\tau(I-A)$. Then
\begin{eqnarray*}
L(\tau(I-A))&=&\mu\circ DT(\tau(I-A))\,=\,\mu\,[\mbox{trace }(I-A\otimes\sum_{m=0}^\infty A^m)]\\
&=&-\mu\,[\mbox{trace }(A\otimes\sum_{m=0}^\infty A^m)]\hspace{2cm}\mbox{since }1\otimes x\mbox{ is a boundary}\\
&=&-\sum_{m=0}^\infty\sum_{i,k=1}^n\mu\,[A_{ik}\otimes A^m_{ki}].
\end{eqnarray*}
It is sufficient to show that
\begin{equation}\label{lside}
\sum_{i,k=1}^n\mu\,[A_{ik}\otimes A^m_{ki}]\,=\,\varepsilon\left(\frac{\mbox{trace }A^{m+1}}{m+1}\right)
\end{equation}
Both sides are clearly $0$ for $\gamma\in\Gamma$ with $\xi(\gamma)\geq 0$. Call the left side of (\ref{lside}) $X$ and let $\gamma\in\Gamma$ with $\xi(\gamma)<0$. Then
\begin{eqnarray*}
X(\gamma)&=&\frac{1}{\xi(\gamma)}\sum_{i_1,\ldots\hspace{-1pt},i_{m+1}=1}^n\sum_{h_1,\ldots\hspace{-1pt},h_{m+1}\in G \atop \{h_1\cdots h_{m+1}\}=\gamma}\xi(h_1)A_{i_1i_2}(h_1)\cdots A_{i_{m+1}i_1}(h_{m+1}).
\end{eqnarray*}
Let $\mathbb{Z}_{m+1}$ act on $\{1,\ldots\hspace{-1pt},m+1\}$ by the cycle $(1\,2\,\cdots\,m+1)$ and on $\{1,\ldots\hspace{-1pt},n\}^{m+1}$ by rotation. For $x\in\{1,\ldots\hspace{-1pt},n\}^{m+1}$ denote by $[x]$ the orbit of $x$ and by $S$ the orbit set. We get
\begin{eqnarray*}
X(\gamma)&=&\frac{1}{\xi(\gamma)}\sum_{i_1,\ldots\hspace{-1pt},i_{m+1}=1}^n\sum_{h_1,\ldots\hspace{-1pt},h_{m+1}\in G \atop \{h_1\cdots h_{m+1}\}=\gamma}\frac{1}{m+1}\sum_{t\in\mathbb{Z}_{m+1}}\xi(h_{t1})A_{i_1i_2}(h_{t1})\cdots A_{i_{m+1}i_1}(h_{t(m+1)})\\
&=&\frac{1}{m+1}\frac{1}{\xi(\gamma)}\sum_{[x]\in S}\sum_{(i_1,\ldots\hspace{-1pt},i_{m+1})\in [x]}\sum_{h_1,\ldots\hspace{-1pt},h_{m+1}\in G \atop \{h_1\cdots h_{m+1}\}=\gamma}\sum_{t\in\mathbb{Z}_{m+1}}\xi(h_{t1})A_{i_1i_2}(h_{t1})\cdots A_{i_{m+1}i_1}(h_{t(m+1)})\\
&=&\frac{1}{m+1}\frac{1}{\xi(\gamma)}\sum_{[x]\in S}\sum_{h_1,\ldots\hspace{-1pt},h_{m+1}\in G \atop \{h_1\cdots h_{m+1}\}=\gamma}\sum_{t\in\mathbb{Z}_{m+1}}\xi(h_{t1})\sum_{(i_1,\ldots\hspace{-1pt},i_{m+1})\in [x]}A_{i_1i_2}(h_{t1})\cdots A_{i_{m+1}i_1}(h_{t(m+1)})
\end{eqnarray*}
Now
\[\sum_{(i_1,\ldots\hspace{-1pt},i_{m+1})\in [x]}A_{i_1i_2}(h_{t1})\cdots A_{i_{m+1}i_1}(h_{t(m+1)})\,=\,\sum_{(i_1,\ldots\hspace{-1pt},i_{m+1})\in [x]}A_{i_1i_2}(h_1)\cdots A_{i_{m+1}i_1}(h_{m+1}),\]
since the orbit is obtained by shifting $(i_1,\ldots\hspace{-1pt},i_{m+1})$, so
\begin{eqnarray*}
X(\gamma)&=&\frac{1}{m+1}\frac{1}{\xi(\gamma)}\sum_{[x]\in S}\sum_{h_1,\ldots\hspace{-1pt},h_{m+1}\in G \atop \{h_1\cdots h_{m+1}\}=\gamma}\sum_{(i_1,\ldots\hspace{-1pt},i_{m+1})\in [x]}A_{i_1i_2}(h_1)\cdots A_{i_{m+1}i_1}(h_{m+1})\sum_{t\in\mathbb{Z}_{m+1}}\xi(h_{t1})\\
&=&\frac{1}{m+1}\frac{\xi(\gamma)}{\xi(\gamma)}\sum_{[x]\in S}\sum_{(i_1,\ldots\hspace{-1pt},i_{m+1})\in [x]}\sum_{h_1,\ldots\hspace{-1pt},h_{m+1}\in G \atop \{h_1\cdots h_{m+1}\}=\gamma}A_{i_1i_2}(h_1)\cdots A_{i_{m+1}i_1}(h_{m+1})\\
&=&\frac{1}{m+1}\sum_{i_1,\ldots\hspace{-1pt},i_{m+1}=1}^n\sum_{h_1,\ldots\hspace{-1pt},h_{m+1}\in G \atop \{h_1\cdots h_{m+1}\}=\gamma}A_{i_1i_2}(h_1)\cdots A_{i_{m+1}i_1}(h_{m+1})\\
&=&\frac{1}{m+1}\,\varepsilon(\mbox{trace }A^{m+1})(\gamma).
\end{eqnarray*}
\end{proof}
In the case where $\xi$ is a homomorphism to the integers it is now easily seen that $\mathfrak{L}$ agrees with Pajitnov's $\mathfrak{L}$
from \cite{pajitn} on $\overline{W}$ once the correct identifications are made.
\section{Geometry of Morse forms}
\subsection{The chain homotopy equivalence}\label{sb4.1}
In \cite{pajirn,pajitn}, Pajitnov defines a condition $(\mathfrak{C}')$, in \cite{pajiov} denoted
by $(\mathfrak{C}\mathcal{C})$, for an $f$-gradient $v$, where $f:M\to S^1$ is a Morse function
that induces a surjection $\xi$ on fundamental group. For the full condition we refer the reader
to these papers, but informally, it can be described as follows: just as in \ref{CHtyTp} we get a
cobordism $(M_N,N,tN)$ and a Morse map $\bar{f}:(M_N,N,tN)\to([0,b],\{0\},\{b\})$.
Here we use an arbitrary $b>0$ instead of just $b=1$ to indicate that $f$ might come from a rational Morse form.
Now the condition $(\mathfrak{C}')$ requires a Morse map $\psi$ on $N$ and a $\psi$-gradient $u$
which gives a handle decomposition on $N$ and $tN$. The vector field $v$ which lifts to a vector field $v'$ on $M_N$ now has to satisfy a ``cellularity condition'':
whenever $p$ is a critical point of $\bar{f}$ of index $i$, it should be the case that some thickening of $D_L(p)$ is attached to the union of the $(i-1)$-handles of $N$ and of $M_N$. 
Also a thickening of an $i$-handle in $tN$ has to flow under $-v'$ into the $i$-skeleton of $N$ and $M_N$. A symmetric condition holds for right hand discs and handles of $N$.\\[0.2cm]
By Pajitnov \cite[Prop.5.4]{pajiov}, the set of $f$-gradients satisfying $(\mathfrak{C'})$ and
the transversality condition is $C^0$-open and dense in the set of $f$-gradients that satisfy the
transversality condition. Such gradients should be thought of as cellular approximations of arbitrary $f$-gradients.\\[0.2cm]
Now let $\rho:\tilde{M}\to M$ be the universal cover, $\tilde{f}:\tilde{M}\to\mathbb{R}$ the lifting
of $f$ and $\tilde{M}_N=\tilde{f}^{-1}([0,b])$. As in \ref{app2} we get a filtration of $M_N$ which
we denote by $M_i$ so that $C_i^{MS}(\tilde{M}_N)=H_i(\tilde{M}_i,\tilde{M}_{i-1})$ gives a free
$\mathbb{Z}H$ complex calculating $H_\ast(\tilde{M}_N)$, where $H=\ker \xi$. We also have
$C_i^{MS}(\tilde{M}_N)=C_i^{MS}(\tilde{N})\oplus C_i^{MS}(\tilde{M}_N,\tilde{N})$. Denote
$C_\ast^{MS}(\tilde{M}_N;G)=\mathbb{Z}G\otimes_{\mathbb{Z}H}C_\ast^{MS}(\tilde{M}_N)$ and similarly
for other $\mathbb{Z}H$ complexes.\\[0.2cm]
The flow of $-v$ on $M_N$ induces a chain map \\$k=\left(\begin{array}{c}k_1\\k_2\end{array}\right):
C_\ast^{MS}(\tilde{N};G)\to C^{MS}_\ast(\tilde{N};G)\oplus C_\ast^{MS}(\tilde{M}_N,\tilde{N};G)=
C_\ast^{MS}(\tilde{M}_N;G)$ by starting in $C_\delta^k(u)\subset t\tilde{N}$ and flowing into
$(\tilde{M}_N)_i$. See appendix \ref{app2} for the sets $C_\delta(\psi)$ and details on this flowing.
Notice that $\delta>0$ is given through condition $(\mathfrak{C}')$. The map $k_1$ is the homological
gradient descent of Pajitnov \cite[\S 4]{pajirn}. Choose a basis of $C_\ast^{MS}(\tilde{N};G)$ by
lifting critical points of $\psi$ and let the matrix $A_i$ represent the homomorphism $k_1$ in
dimension $i$. Since we can choose the liftings within $\tilde{N}$, we get $\|A_{jk}\|<1$ for the
entries of every matrix $A_i$.
\begin{prop}\label{ptor}
Let $v$ be an $f$-gradient satisfying the transversality condition and $(\mathfrak{C}')$. Then
$\varphi(v)$ is a chain homotopy equivalence and
\[\tau(\varphi(v))=\sum_{i=0}^{n-1}(-1)^{i+1}\tau(I-A_i)\in\overline{K}_1^G(\widehat{\mathbb{Z}G}_\xi).\]
\end{prop}
\begin{proof}
Let $\Delta$ be a triangulation of $M$ which has $N$ as a subcomplex
$\Delta'$. This induces a triangulation $\Delta^c$ of $M_N$ which has two copies of $\Delta'$ as
subcomplexes. Denote the one corresponding to $N$ by $\Delta_0$ and the one corresponding to $tN$
by $\Delta_1$. Assume $\Delta$ has the following properties:
\begin{enumerate}
\item $\Delta$ is adjusted to $v$.
\item $\Delta'$ is adjusted to $u$.
\item There is an $\varepsilon>\delta$ such that if $\sigma$ is a $k$-simplex in $\Delta'$,
then $\sigma\subset C^k_\varepsilon(u)$.
\item There is an $\varepsilon>\delta$ such that if $x\in M_N^{(k)}$ and the trajectory of $-v$ starting
at $x$ ends in $fl(x)\in N$, then $fl(x)\in C_\varepsilon^k(u)$.
\end{enumerate}
The existence of such a triangulation is shown in \ref{lastap}.\\[0.2cm]
\ref{app2} gives a simple chain homotopy equivalence
\[\varphi:C^\Delta_\ast(\tilde{M}_N;G)\to C^{MS}_\ast(\tilde{M}_N;G).\]
Let us define a chain map $s:C_\ast^\Delta(\tilde{N};G)
\to C_\ast^\Delta(\tilde{M}_N,\tilde{N};G)\subset C_\ast^\Delta(\tilde{M}_N;G)$. If $\sigma$ is
a simplex in $\Delta'$, look at the liftings $\sigma_0\subset\tilde{N}$ and $\sigma_1\subset t\tilde{N}$
used for the basis of $C_\ast^\Delta(\tilde{M}_N)$. There is exactly one $g\in G$ such that $g\sigma_1=
\sigma_0$ in $\tilde{M}$. Set $s(\sigma_0)=g\sigma_1$. Look at the diagram
\[\begin{array}{ccc}
C^\Delta_\ast(\tilde{N};G)&\stackrel{s}{\longrightarrow}&C^\Delta_\ast(\tilde{M}_N;G)\\
\Big\downarrow\varphi_1& &\Big\downarrow\varphi\\
C^{MS}_\ast(\tilde{N};G)&\stackrel{k}{\longrightarrow}&C_\ast^{MS}(\tilde{M}_N;G)\end{array}\]
Because of property 4.\ above, this diagram commutes. Therefore the map $\left(\begin{array}{cc}
\varphi&0\\0&\varphi_1\end{array}\right)$ is a simple homotopy equivalence between the mapping cones
$\mathcal{C}(i-s)$ and $\mathcal{C}(i-k)$, where $i$ represents inclusion. But by the Deformation
Lemma of Farber and Ranicki \cite{farran}, $\mathcal{C}(i-s)$ is chain homotopy equivalent to
${\rm coker}(i-s)$, in fact simple homotopy equivalent by Ranicki \cite[Prop.1.9]{ranick} (the
corresponding matrix term is just $I$). But ${\rm coker}(i-s)$ is easily seen to be $C_\ast^\Delta
(\tilde{M})$.\\[0.2cm]
After tensoring with the Novikov ring we have the following sequence of chain homotopy equivalences
\[\widehat{\mathbb{Z}G}_\xi\otimes_{\mathbb{Z}G}C_\ast^\Delta(\tilde{M})\to
\widehat{\mathbb{Z}G}_\xi\otimes_{\mathbb{Z}G}\mathcal{C}(i-s)\to
\widehat{\mathbb{Z}G}_\xi\otimes_{\mathbb{Z}G}\mathcal{C}(i-k)\to
{\rm coker}({\rm id}_{\widehat{\mathbb{Z}G}_\xi}\otimes i-k)\]
and all except the last one are simple. The first map is described in the proof of Ranicki
\cite[Prop.1.7]{ranick}. Because of the special form of the vector field $v$ the Novikov complex
$C_\ast(f,v)$ can be identified with ${\rm coker}({\rm id}_{\widehat{\mathbb{Z}G}_\xi}\otimes i-k)$,
see Ranicki \cite[Remark 4.8]{ranick} and Pajitnov \cite[Remark 7.3]{pajirn}.
We claim that this composition is exactly $\varphi(v)$. Denote
the composition by $\theta$.\\[0.2cm]
We denote $t^k\tilde{M}_N=\tilde{f}^{-1}[bk,b(k+1)]$.
Let $\sigma\in C^\Delta_k(\tilde{M})$, lift it to $\bar{\sigma}\in C^\Delta_k(\tilde{M}_N,\tilde{N}
;G)$ (if $\sigma$ is a cell in $\tilde{N}$, lift it to $\bar{\sigma}\subset t\tilde{N}$). Then
\[\theta(\sigma)=[\varphi(0,\bar{\sigma})]=[\varphi_N(\bar{\sigma}),\varphi_r(\bar{\sigma})]\in
{\rm coker}(i-k)=C_\ast(\omega,v),\]
where $\varphi_N(\bar{\sigma})\in C^\Delta_k(\tilde{N};G)$ and $\varphi_r(\bar{\sigma})\in
C^\Delta_k(\tilde{M}_N,\tilde{N};G)$ are defined to give \\$(\varphi_N(\bar{\sigma}),
\varphi_r(\bar{\sigma}))=\varphi(0,\bar{\sigma})$.\\[0.2cm]
Now $[\varphi_N(\bar{\sigma}),\varphi_r(\bar{\sigma})]=[\varphi_N(\bar{\sigma}),0]+[0,\varphi_r(
\bar{\sigma})]$. $\varphi_r(\bar{\sigma})$ represents the part of $\bar{\sigma}$ that flows into
critical points of index $k$ in $\tilde{M}_N$ under $-v$ while $\varphi_N(\bar{\sigma})$ represents
the part that flows into $\tilde{N}$. Now
\begin{eqnarray*}
(i-k)(\sum_{m=0}^\infty k_1^m(x))&=&(\sum_{m=0}^\infty k_1^m(x),0)-
(\sum_{m=0}^\infty k_1^{m+1}(x),k_2(\sum_{m=0}^\infty k_1^m(x)))\\
&=&(x,k_2(\sum_{m=0}^\infty k_1^m(x))).
\end{eqnarray*}
With $x=\varphi_N(\bar{\sigma})$ we thus get
\[[\varphi_N(\bar{\sigma}),0]=[0,k_2(\sum_{m=0}^\infty k_1^m(\varphi_N(\bar{\sigma}))].\]
But $k_2(k_1^m(\varphi_N(\bar{\sigma})))$ represents the part of $\bar{\sigma}$ that flows into
critical points of index $k$ in $t^{-m-1}\tilde{M}_N$ under $-v$. Therefore $\theta=\varphi(v)$
and $\varphi(v)$ is a chain homotopy equivalence whose torsion is given by (\ref{torsion}).
\end{proof}
\begin{rema}\em
Pajitnov \cite{pajirn,pajitn} obtains a chain homotopy equivalence \\$\psi(v):C_\ast(f,v)\to
\widehat{\mathbb{Z}G}_\xi\otimes_{\mathbb{Z}G}C_\ast^\Delta(\tilde{M})$ by including the Novikov
complex into a complex $C'$ simple homotopy equivalent to $\widehat{\mathbb{Z}G}_\xi\otimes_{
\mathbb{Z}G}\mathcal{C}(i-k)$, in fact the map from $C'$ to $\widehat{\mathbb{Z}G}_\xi\otimes_{\mathbb{Z}G}
\mathcal{C}(i-k)$ is just a simple change of basis, compare \cite[\S 7.4]{pajirn}. The composition of this
equivalence with Ranicki's equivalence $\widehat{\mathbb{Z}G}_\xi\otimes_{\mathbb{Z}G}\mathcal{C}
(i-k)\to C_\ast(f,v)$ is readily seen to be the identity on $C_\ast(f,v)$.
\end{rema}
\subsection{Approximation of irrational forms by rational forms}\label{sb4.2}
In this subsection we describe a useful method due to Pajitnov \cite[\S 2B]{pajisp}.
Given a Morse form $\omega$ and an $\omega$-gradient $v$, the induced homomorphism $\bar{\xi}_\omega:H_1(M)\to\mathbb{R}$ splits $H_1(M)$ as $\mathbb{Z}^k\oplus\ker\bar{\xi}_\omega$. Choose $g_1,\ldots\hspace{-1pt},g_k\in G$ so that the images $\bar{g}_1,\ldots\hspace{-1pt},\bar{g}_k\in H_1(M)$ generate the $\mathbb{Z}^k$ part. Now let $\omega_1,\ldots\hspace{-1pt},\omega_k$ be closed $1$-forms with $\bar{\xi}_{\omega_j}(\bar{g}_i)=\delta_{ji}$ and $\bar{\xi}_{\omega_j}$ vanishes on $\ker \bar{\xi}_\omega$.
Then $\xi_{\omega_j}:G\to\mathbb{Z}$ vanishes on $\ker\xi_\omega$ and satisfies $\xi_{\omega_j}(g_i)=\delta_{ji}$. Furthermore the closed $1$-forms can be chosen to vanish in a neighborhood of the critical points of $\omega$.\\[0.2cm]
For $x\in\mathbb{R}^k$ we can now define $\omega_x=\omega+\sum\limits_{j=1}^k x_j\omega_j$.
By choosing the $x_j$ small we can make sure that the $\omega$-gradient $v$ is also an $\omega_x$-gradient. To see this notice that in a neighborhood of the critical points of $\omega$ the new form agrees with $\omega$.
Denote the complement of this neighborhood by $C$.
Because of the compactness of $C$ and Lemma \ref{smalo} there is a $K>0$ such that $\omega_p(v(p))\geq K$ for all $p\in C$. 
Now the $x_j$ have to be chosen so small that $(\omega_x)_p(v(p))>0$ for all $p\in C$ which is possible again by compactness. Lemma \ref{smalo} now gives that $v$ is an $\omega_x$-gradient.\\[0.2cm]
We have $\xi_{\omega_x}(g_j)=\xi_\omega(g_j)+x_j\xi_{\omega_j}(g_j)$, so we can also choose the $x_j$ to have $\xi_{\omega_x}:G\to\mathbb{R}$ factor through $\mathbb{Q}$. Hence we get
\begin{lem}\label{approx}For a Morse form $\omega$ and an $\omega$-gradient $v$ there exists a rational Morse form $\omega'$ with the same set of critical points and that agrees with $\omega$ in a neighborhood of these critical points such that $v$ is also an $\omega'$-gradient.
\end{lem}
Let us compare the Novikov complexes we obtain for a Morse form $\omega$ and a rational approximation $\omega'$ that both use the same vector field $v$. The complexes are taken over different rings, $\widehat{\mathbb{Z}G}_{\xi_\omega}$ and $\widehat{\mathbb{Z}G}_{\xi_{\omega'}}$ respectively.
But for two critical points $p,q$ of adjacent index the elements $\tilde{\partial}(p,q)\in\widehat{\mathbb{Z}G}_{\xi_\omega}$ and $\tilde{\partial}'(p,q)\in\widehat{\mathbb{Z}G}_{\xi_{\omega'}}$ agree when viewed as elements of $\widehat{\widehat{\mathbb{Z}G}}$ since both count the number of flowlines between $\tilde{p}$ and translates of $\tilde{q}$, and these only depend on $v$.
So we can compare chain complexes even though they are over different rings. This is an important observation and will remain useful in the next subsection.
\subsection{Comparison of the eta function with torsion}
Again let $\omega$ be a Morse form. An $\omega$-gradient $v$ satisfies the condition
$(\mathfrak{AC})$, if there exists a rational Morse form $\omega'$ such that $v$ is an
$\omega'$-gradient and as such it satisfies $(\mathfrak{C}')$. We think of this condition
as ``approximately $(\mathfrak{C}')$''.\\[0.2cm]
We want to carry over the density results of Pajitnov \cite{pajiov}.
Then $C^0$-density in $\mathcal{G}(\omega)$ can be seen as follows:
given an $\omega$-gradient $v'$ there is by Lemma \ref{approx} a rational Morse form $\omega'$ that agrees with $\omega$ near the critical points and such that $v'$ is also an $\omega'$-gradient.
Now the density of $\omega'$-gradients  satisfying $(\mathfrak{C}')$ allows us to choose a vector field $v$ as close as we like to $v'$. To see that we can find an $\omega$-gradient this way observe that 1. of Lemma \ref{smalo} is trivially fulfilled and since $\omega(v')\geq K>0$ away from a neighborhood of the critical points we get $\omega(v)>0$ by choosing $v$ close enough to $v'$.
Therefore $v$ is an $\omega$ gradient satisfying $(\mathfrak{AC})$. The $C^0$-openness now
follows from Pajitnov \cite[Prop.5.4]{pajiov}.\\[0.2cm]
Now if an $\omega$-gradient $v$ satisfies $(\mathfrak{AC})$, let $\omega'$ be the rational Morse
form as in the definition and denote by $\xi:G\to\mathbb{R}$ and $\xi':G\to\mathbb{R}$ the
homomorphisms induced by $\omega$ and $\omega'$. For the rational form $\omega'$ we can form
$\tilde{N}$, the $\mathbb{Z}G$ complex $\mathcal{C}(i-k)$, the homomorphism $k_1:C_i^{MS}(
\tilde{N};G)\to C_i^{MS}(\tilde{N};G)$ and the matrix $A_i$ just as in \ref{sb4.1} and the proof
of Proposition \ref{ptor}.
If we can show that $I-A_i$ is invertible over $\widehat{\mathbb{Z}G}_\xi$ we get the chain homotopy
equivalence between $\widehat{\mathbb{Z}G}_\xi\otimes_{\mathbb{Z}G} \mathcal{C}(i-k)$ and coker$($id$\otimes i-k)$.
We know from 4.1 that $I-A_i$ is invertible over $\widehat{\mathbb{Z}G}_{\xi'}$ and that the
cokernel over this Novikov ring is exactly the Novikov complex $C_\ast(\omega',v)$.
In order to keep the notation simple denote the matrix $A_i$ by $B$.
\begin{rema}\label{apprem}\em
That $I-B$ is invertible over $\widehat{\mathbb{Z}G}_{\xi'}$ is easily seen since a basis can be
chosen so that $\|B_{ij}\|_{\xi'}<1$ for every entry $B_{ij}$. If we could choose a basis of $D$
such that $\|B_{ij}\|_\xi<1$ we would also immediately get that $I-B$ is invertible over
$\widehat{\mathbb{Z}G}_\xi$. A similar argument is used in Latour \cite[\S 2.23]{latour}. But since
we have to choose liftings of cells instead of critical points it is not clear that a nice basis
can be chosen. So instead of trying to find a nice basis we use the notion of $\xi$-regular
matrices.
\end{rema}
\begin{prop}The matrix $B$ over $\mathbb{Z}G$ is $\xi$-regular.\label{regul}
\end{prop}
\begin{proof}
We have chosen a basis of $D_i$ by choosing $i$-cells in $\tilde{N}$, call these cells $\sigma_k$. If $h\in$ supp $B_{jk}$, then there exist negative flowlines from $\sigma_j$ to $h\sigma_k$ by the construction of $B$.\\[0.2cm]
We need to show that there exists a $K<0$ with the property that given $m\geq 1$ and indices $j,n_1,\ldots\hspace{-1pt},n_{m-1}$ and $g_1\in$ supp $B_{jn_1}$, $g_2\in$ supp $B_{n_1n_2},\ldots\hspace{-1pt},g_m\in$ supp $B_{n_{m-1}j}$ we have $\xi(g_1\cdots g_m)\leq K\cdot m$.\\[0.2cm]
Now we have to recall the proof of the Main Theorem in Pajitnov \cite[\S 5]{pajitn}.
Every cell $\sigma_k$ defines a thickened sphere in $\coprod\limits_{l\in\mathbb{Z}}\tilde{V}_l^{[i]}/\tilde{V}_l^{(i-1)}$ that we denote by $s_k$, also let $g=g_1\cdots g_m$.
The spaces $\tilde{V}_l^{[i]}$ and $\tilde{V}_l^{(i)}$ are defined in Pajitnov \cite[\S 4.4,\S 4.5]{pajitn}.
Since $g_1\in$ supp $B_{jn_1}$ we have that $-\tilde{v}$ induces a homologically nontrivial map from $s_j$ to $g_1s_{n_1}$. Similarly every $g_l\in$ supp $B_{n_{l-1}n_l}$ gives rise to a homologically nontrivial map from $g_1\cdots g_{l-1}s_{n_{l-1}}$ to $g_1\cdots g_ls_{n_l}$. 
The composition of all these maps plus $g^{-1}:gs_j\to s_j$ is homologically nontrivial and hence has a fixed point other than the base point, compare the proof of Lemma 5.1 in \cite{pajitn}.
Notice that the existence of this fixed point does not require the condition that closed orbits of $v$ are nondegenerate.
This fixed point corresponds to a flow line $\gamma:[a_1,a_2]\to\tilde{M}$ of $-\tilde{v}$ with $\gamma(a_1)=x\in\sigma_j$ and $\gamma(a_2)=gx\in g\sigma_j$ which passes through the cells $g_1\cdots g_l\sigma_{n_l}$.\\[0.2cm]
We need the following 
\begin{lem}There exists a $K<0$ such that for every flowline $\gamma$ of $-\tilde{v}$ that starts in $\tilde{N}_0$ and ends in $\tilde{N}_{-1}$ we have $\int_{\rho\circ\gamma}\omega\leq K$.
\label{easyl}
\end{lem}
\begin{proof}
We have $\rho^\ast\omega'=d\tilde{f}'$ and $\tilde{N}_k=(\tilde{f}')^{-1}(k\cdot b)$ with $b$ as in 4.1. 
Since $\omega'$ is rational and $\tilde{f}'$ has no critical points in $\tilde{N}_0$ there is a $t<0$ such that $(\tilde{f}')^{-1}([t,0])$ also contains no critical points.
So if $\gamma_p$ is a flowline of $-\tilde{v}$ with $\gamma_p(0)=p\in\tilde{N}_0$, there is a $t_p>0$ which depends smoothly on $p$ such that $\gamma_p(t_p)\in(\tilde{f}')^{-1}(\{t\})$. Now
\[\int\limits_{\rho\circ\gamma_p|[0,t_p]}\omega=\int\limits_0^{t_p}\omega_{\rho\circ\gamma_p(s)}(-v(\rho\circ\gamma_p(s)))=-\int\limits_0^{t_p}\omega_{\rho\circ\gamma_p(s)}(v(\rho\circ\gamma_p(s)))<0\]
by Lemma \ref{smalo}. Since $\ker\xi'$ acts cocompactly on $\tilde{N}_0$ and the value of the integral depends smoothly on $p\in\tilde{N}_0$ there is $K<0$ such that 
\[\int\limits_{\rho\circ\gamma_p|[0,t_p]}\omega\leq K.\]
 This $K$ now works for the Lemma, since integrating over a longer flowline will just make the integral smaller.
\end{proof}
Conclusion of the proof of \ref{regul}:
our flowline $\gamma$ is the concatenation of flowlines $\gamma_1,\ldots\hspace{-1pt},\gamma_m$ to which Lemma \ref{easyl} applies. Let $\tilde{f}:\tilde{M}\to\mathbb{R}$ satisfy $\rho^\ast\omega=d\tilde{f}$. Then we get
\[\xi(g)=\tilde{f}(gx)-\tilde{f}(x)=\int\limits_\gamma d\tilde{f}=\sum_{l=1}^m\int\limits_{\gamma_l}d\tilde{f}=\sum_{l=1}^m\int\limits_{\gamma_l}\rho^\ast\omega=\sum_{l=1}^m\int\limits_{\rho\circ\gamma_l}\omega\leq m\cdot K.\]
Therefore $\xi(g_1)+\ldots+\xi(g_m)\leq m\cdot K$ for all $g_1\in$ supp $B_{jn_1},\ldots\hspace{-1pt},g_m\in$ supp $B_{n_{m-1}j}$ which implies that $B$ is $\xi$-regular.
\end{proof}
Define $\mathcal{G}_0(\omega)=\{v\in\mathcal{G}(\omega)\,|\,v$ satisfies $(\mathfrak{AC})\}$. By the remarks above, $\mathcal{G}_0(\omega)$ is $C^0$-dense in $\mathcal{G}(\omega)$.
\begin{theo}\label{maint}
Let $\omega$ be a Morse form on a smooth connected closed manifold $M^n$. Let $\xi:G\to\mathbb{R}$
be induced by $\omega$ and let $C^\Delta_\ast(\tilde{M})$ be the simplicial $\mathbb{Z}G$ complex
coming from a smooth triangulation of $M$. For every $v\in\mathcal{G}_0(\omega)$ there is a natural
chain homotopy equivalence $\varphi(v):\widehat{\mathbb{Z}G}_\xi\otimes_{\mathbb{Z}G}C_\ast^\Delta
(\tilde{M})\to C_\ast(\omega,v)$ given by (\ref{formula}) whose torsion $\tau(\varphi(v))$ lies in
$\overline{W}$ and satisfies
\[\mathfrak{L}(\tau(\varphi(v)))=\eta(-v).\]
\end{theo}
\begin{proof}
Since $v\in\mathcal{G}_0(\omega)$ we can form the $\mathbb{Z}G$ complex $\mathcal{C}(i-k)$ from
the proof of Proposition \ref{ptor} which is simple homotopy equivalent to $C^\Delta_\ast(\tilde{M})$.
The matrices $A_i$ are $\xi$-regular by Proposition \ref{regul}, so the projection of
$\widehat{\mathbb{Z}G}_\xi\otimes_{\mathbb{Z}G}\mathcal{C}(i-k)\to$ coker$($id$\otimes i-k,\widehat{\mathbb{Z}G}_\xi)$ is a chain homotopy equivalence.
We have already seen that the boundary homomorphisms of $C_\ast(\omega,v)$ and $C_\ast(\omega',v)$ are the same when viewed as matrices over $\widehat{\widehat{\mathbb{Z}G}}$.
The same holds for coker$($id$\otimes i-k,\widehat{\mathbb{Z}G}_\xi)$ and coker$($id$\otimes i-k,\widehat{\mathbb{Z}G}_{\xi'})$. But since we identified coker$($id$\otimes i-k,\widehat{\mathbb{Z}G}_{\xi'})$ with $C_\ast(\omega',v)$ we now get that coker$($id$\otimes i-k,\widehat{\mathbb{Z}G}_\xi)$ is the same complex as $C_\ast(\omega,v)$.
Also, if we use the triangulation from the proof of Proposition \ref{ptor} we get that $\varphi(v)$
factors as $p\circ s$, where $s$ is simple. Therefore
\[\tau(\varphi(v))=\sum_{i=0}^{n-1}(-1)^{i+1}\tau(I-A_i)\in\overline{K}^G_1(\widehat{\mathbb{Z}G}_\xi).\]
By Proposition \ref{logar} we have 
\[\mathfrak{L}(\tau(\varphi(v)))=\sum_{i=0}^{n-1}(-1)^i\sum_{m=1}^\infty\frac{\varepsilon(\mbox{trace }A^m)}{m}.\]
By the proof of the Main Theorem in \S 5 of \cite{pajitn} the right hand side is exactly $\eta(-v)$. Of course, \cite{pajitn} only shows this in the rational case, but $\eta(-v)$ is independent of the homomorphism $\xi$ when viewed as an element of $\widehat{\widehat{\mathbb{Q}\Gamma}}$.
\end{proof}
\section{Comparison with Reidemeister Torsion}
As mentioned in the introduction, for singular vector fields one of the first formulas to relate the torsion of the Novikov complex to a zeta function appeared in Hutchings and Lee \cite{hutlee,hutle2}, a generalization appeared in Hutchings \cite{hutcth,hutchi} looking similar to Theorem \ref{maint}, but using quite different methods.
In this section we will relate these results, in fact we will show that Theorem \ref{maint} implies \cite[Theorem B]{hutchi}, at least for gradients satisfying condition $(\mathfrak{AC})$.\\[0.2cm]
All these papers deal with commutative invariants only, so let $\overline{M}$ be the universal abelian cover of $M$ and $H=H_1(M)$ the covering transformation group.
Let $\omega$ be a Morse form and $v$ an $\omega$-gradient satisfying the transversality condition. The Novikov complex in Hutchings \cite{hutchi} is given by 
$\overline{C}_\ast(\omega,v)=\widehat{\mathbb{Z}H}_{\bar{\xi}}\otimes_{\widehat{\mathbb{Z}G}_\xi}C_\ast(\omega,v)$. 
Similarly $\widehat{\mathbb{Z}H}_{\bar{\xi}}\otimes_{\mathbb{Z}H}C^\Delta_\ast(\overline{M})=\widehat{\mathbb{Z}H}_{\bar{\xi}}\otimes_{\mathbb{Z}G} C^\Delta_\ast(\tilde{M})$ and therefore a chain equivalence
$\varphi(v):\widehat{\mathbb{Z}G}_\xi\otimes_{\mathbb{Z}G}C^\Delta_\ast(\tilde{M})\to C_\ast(\omega,v)$ induces a chain equivalence $\bar{\varphi}(v)=$ id$\,\otimes_{\widehat{\mathbb{Z}G}_\xi}\varphi(v):\widehat{\mathbb{Z}H}_{\bar{\xi}}\otimes_{\mathbb{Z}H}C^\Delta_\ast(\overline{M})\to\overline{C}_\ast(\omega,v)$ with
$\tau(\bar{\varphi}(v))=\varepsilon_\ast\tau(\varphi(v))\in\overline{K}^H_1(\widehat{\mathbb{Z}H}_{\bar{\xi}})$.\\[0.2cm]
Let $Q\widehat{\mathbb{Z}H}_{\bar{\xi}}$ be the localization of $\widehat{\mathbb{Z}H}_{\bar{\xi}}$ along non-zero divisors. It is known that $Q\widehat{\mathbb{Z}H}_{\bar{\xi}}$ is a finite direct product of fields, $Q\widehat{\mathbb{Z}H}_{\bar{\xi}}=\bigoplus\limits_{j=1}^k F_j$, see Hutchings \cite[Lemma A.4]{hutchi} or Geoghegan and Nicas \cite[Lemma 7.8]{geonic}. Denote $p_j:Q\widehat{\mathbb{Z}H}_{\bar{\xi}}\to F_j$ for the projection.\\[0.2cm]
To define torsion in the sense of Hutchings \cite{hutchi}, we need one more construction. Let $R$ be a commutative ring with unit and $U$ a subgroup of $R^\ast$, the group of units of $R$. We say that two elements of $R$ are \em equivalent\em, $r\sim s$, if there exists a $u\in U$ such that $ru=s$. We denote by $R/U$ the set of equivalence classes. The multiplication on $R$ turns $R/U$ into a semigroup which contains $R^\ast/U$ as a subgroup.
\begin{defi}\em \cite[Def.A.1]{hutchi} Let $F$ be a field and $C_\ast$ a finite complex over $F$ with a fixed basis and $U$ a subgroup of $F^\ast$. Then the \em Reidemeister torsion \em of $C_\ast$ is defined to be
\[\tau_R(C_\ast,U)=\left\{\begin{array}{rll}0&\in F/U&\mbox{if }C_\ast\mbox{ is not acyclic}\\ \det(\tau(C_\ast))^{-1}&\in F^\ast/U\subset F/U&\mbox{if }C_\ast\mbox{ is acyclic}\end{array}\right.\]
Here $\tau(C_\ast)\in K_1(F)$ is Whitehead torsion.
\end{defi}
We take the inverse of the determinant, because Hutchings \cite{hutchi} uses a different sign convention for torsion. When the group of units is clear, we will suppress it in the notation
of the torsion.\\[0.2cm]
Now $\pm H$ is a subgroup of $Q\widehat{\mathbb{Z}H}_{\bar{\xi}}^\ast$. Denote $H_j=p_j(H)\subset F_j^\ast$ for $j=1,\ldots\hspace{-1pt},k$. Then \\$Q\widehat{\mathbb{Z}H}_{\bar{\xi}}/\pm H=\bigoplus\limits_{j=1}^k F_j/\pm H_j$. Hutchings \cite[\S 1.5]{hutchi} defines
\[T_m=\sum_{j=1}^k\tau_R(F_j\otimes_{\widehat{\mathbb{Z}H}_{\bar{\xi}}}\overline{C}_\ast(\omega,v))\in Q\widehat{\mathbb{Z}H}_{\bar{\xi}}/\pm H\]
and 
\[T(\overline{M})=\sum_{j=1}^k\tau_R(F_j\otimes_{\mathbb{Z}H}C^\Delta_\ast(\overline{M}))\in Q\widehat{\mathbb{Z}H}_{\bar{\xi}}/\pm H.\]
Notice that $T(\overline{M})$, and hence $T_m$, can only be nonzero if $\chi(M)=0$.\footnote{The
Euler characteristic of the complex $F_j\otimes_{\mathbb{Z}H}C^\Delta_\ast(\overline{M})$ equals
the Euler characteristic of $M$ and since $F_j$ is a field it can be calculated from the homology of that complex.}
These torsions are related by the zeta function of $-v$ which is defined as follows. Let \\$\widehat{RG}^-_\xi=\{\lambda\in\widehat{RG}_\xi\,|\,\|\lambda\|<1\}$. Similarly we get $\widehat{RH}^-_{\bar{\xi}}$ and $\widehat{R\Gamma}^-_\xi$.
Define $\log:1+\widehat{\mathbb{Q}H}^-_{\bar{\xi}}\to\widehat{\mathbb{Q}H}^-_{\bar{\xi}}$ and $\exp:\widehat{\mathbb{Q}H}^-_{\bar{\xi}}\to 1+\widehat{\mathbb{Q}H}^-_{\bar{\xi}}$ by
\[\log(1+a)=\sum_{m=1}^\infty (-1)^{m+1}\frac{a^m}{m}\hspace{1cm}\mbox{and}\hspace{1cm}\exp(a)=\sum_{m=0}^\infty\frac{a^m}{m!}.\]
It is readily seen that $\log$ and $\exp$ are well defined and mutually inverse to each other.\\[0.2cm]
If $v\in\mathcal{G}(\omega)$, then $\eta(-v)$ is an element of $\widehat{\mathbb{Q}\Gamma}^-_\xi$ and we define the \em zeta function \em of $-v$ to be 
\[\zeta(-v)=\exp(\varepsilon(\eta(-v)))\in 1+\widehat{\mathbb{Q}H}^-_{\bar{\xi}}.\]
Hutchings relates $T_m$ and $T(\overline{M})$ by
\begin{theo}\cite[Theorem B]{hutchi}\label{htheo} Let $\omega$ be a Morse form on the closed connected smooth manifold $M$, let $v\in\mathcal{G}(\omega)$ and let $\iota:1+\widehat{\mathbb{Q}H}^-_{\bar{\xi}}\to Q\widehat{\mathbb{Z}H}_{\bar{\xi}}/\pm H$ be given by inclusion and projection. Then $T_m\cdot\iota(\zeta(-v))=T(\overline{M})$.
\end{theo}
We will show how to derive Theorem \ref{htheo} from Theorem \ref{maint} for $v\in\mathcal{G}_0(\omega)$. It would be desirable to extend Theorem \ref{maint} to $v\in\mathcal{G}(\omega)$, possibly by the methods of \cite{hutchi}.\\[0.2cm]
So let $\varphi(v):\widehat{\mathbb{Z}G}_\xi\otimes_{\mathbb{Z}G}C_\ast^\Delta(\tilde{M})\to C_\ast(\omega,v)$ be the chain homotopy equivalence from Theorem \ref{maint}. As seen above, this induces a chain homotopy equivalence \\$\bar{\varphi}(v):\widehat{\mathbb{Z}H}_\xi\otimes_{\mathbb{Z}H}C_\ast^\Delta(\overline{M})\to \overline{C}_\ast(\omega,v)$.
The homomorphism $\mathfrak{L}$ actually gives a map $\mathfrak{L}:\overline{W}\to\widehat{\mathbb{Q}\Gamma}^-_\xi$. It is easy to see that the following diagram commutes.
\[\begin{array}{ccccccc}
\overline{W}&\stackrel{\mathfrak{L}}{\longrightarrow}&\widehat{\mathbb{Q}\Gamma}^-_\xi&\stackrel{\varepsilon}{\longrightarrow}&\widehat{\mathbb{Q}H}^-_{\bar{\xi}}&\stackrel{\exp}{\longrightarrow}&1+\widehat{\mathbb{Q}H}^-_{\bar{\xi}}\\[0.3cm]
\Big\downarrow&\multicolumn{5}{c}{ }&\Big\downarrow\\[0.3cm]
\overline{K}_1^G(\widehat{\mathbb{Z}G}_\xi)&\stackrel{\varepsilon_\ast}{\longrightarrow}&\multicolumn{3}{c}{\overline{K}_1^H(\widehat{\mathbb{Z}H}_{\bar{\xi}})}&\stackrel{\det}{\longrightarrow}&\widehat{\mathbb{Q}H}^\ast_{\bar{\xi}}/\pm H\end{array}\]
The last vertical arrow is inclusion of $1+\widehat{\mathbb{Q}H}^-_{\bar{\xi}}$ into $\widehat{\mathbb{Q}H}^\ast_{\bar{\xi}}$ followed by projection.
Hence we get $\det(\tau(\bar{\varphi}(v)))=\bar{\zeta}(-v)\in\widehat{\mathbb{Q}H}^\ast_{\bar{\xi}}/\pm H$. Look at the commutative diagram
\begin{equation}
\label{commd}
\begin{array}{ccccc}
\overline{K}_1^H(\widehat{\mathbb{Z}H}_{\bar{\xi}})&\stackrel{\iota_\ast}{\longrightarrow}&\overline{K}_1^H(Q\widehat{\mathbb{Z}H}_{\bar{\xi}})&\stackrel{\simeq}{\longrightarrow}&\bigoplus\limits_{j=1}^k\overline{K}_1^{H_j}(F_j)\\[0.4cm]
\Big\downarrow\det& &\Big\downarrow\det& &\Big\downarrow\oplus\det\\[0.3cm]
\widehat{\mathbb{Q}H}^\ast_{\bar{\xi}}/\pm H&\stackrel{\bar{\iota}}{\longrightarrow}&Q\widehat{\mathbb{Z}H}^\ast_{\bar{\xi}}/\pm H&\stackrel{\simeq}{\longrightarrow}&\bigoplus\limits_{j=1}^k F_j^\ast/\pm H_j
\end{array}\end{equation}
We will show that \begin{equation}p_j(T_m\cdot\iota(\zeta(-v)))=p_j(T(\overline{M}))\hspace{0.4cm}\mbox{for every }j=1,\ldots\hspace{-1pt},k.\label{formu}\end{equation}
We have to compare $\det\circ (p_j)_\ast\circ\iota_\ast\tau(\bar{\varphi}(v))$ with $\tau_R(F_j\otimes\overline{C}(\omega,v))$ and $\tau_R(F_j\otimes C^\Delta_\ast(\overline{M}))$.\\[0.2cm]
If $F_j\otimes\overline{C}(\omega,v)$ is acyclic, then so is $F_j\otimes C_\ast^\Delta(\overline{M})$, because id$_{F_j}\otimes\bar{\varphi}(v)$ is an equivalence of these complexes. Furthermore
\[\tau(\mbox{id}_{F_j}\otimes\bar{\varphi}(v))=\tau(F_j\otimes\overline{C}(\omega,v))-\tau(F_j\otimes C^\Delta_\ast(\overline{M}))\in \overline{K}_1^{H_j}(F_j).\]
So $\det(\tau($id$_{F_j}\otimes\bar{\varphi}(v)))=\tau_R(F_j\otimes C_\ast^\Delta(\overline{M}))\cdot\tau_R(F_j\otimes\overline{C}(\omega,v))^{-1}$. This gives by (\ref{commd})
and Theorem \ref{maint}
\[\tau_R(F_j\otimes\overline{C}(\omega,v))\cdot p_j(\zeta(-v))=\tau_R(F_j\otimes C_\ast^\Delta(\overline{M})).\]
If  $F_j\otimes\overline{C}(\omega,v)$ is not acyclic, neither is $F_j\otimes C_\ast^\Delta(\overline{M})$ and (\ref{formu}) reduces to $0=0$. Hence we obtain the desired formula $T_m\cdot\iota(\zeta(-v))=T(\overline{M})$.\\[0.2cm]
Since the complexes $F_j\otimes\overline{C}(\omega,v)$ do not always have to be acyclic, Theorem \ref{htheo} cannot recover the zeta function in general just from the torsion information.
In particular, Theorem \ref{htheo} contains no information for $\chi(M)\not=0$. To see that we can
get reasonable results for $\chi(M)\not=0$ we have the following
\begin{exam}\em Let $M$ be a closed surface of genus 2 and $f:M\to S^1$ a Morse map indicated by Figure \ref{first},
\begin{figure}[ht]
\begin{center}
\epsfig{file=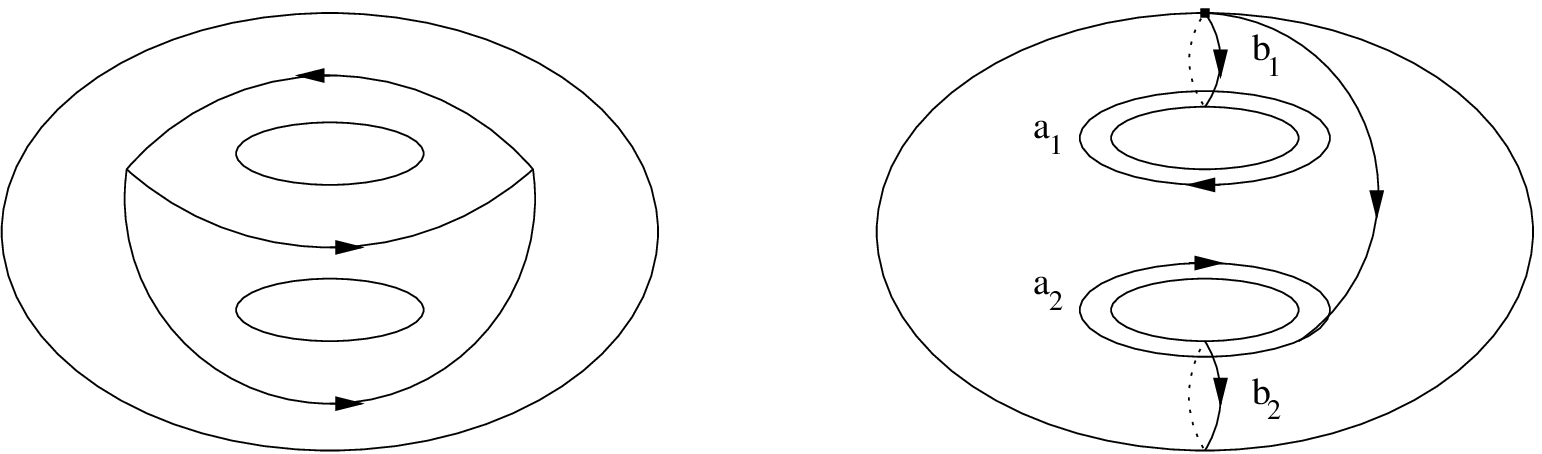,height=4cm}
\end{center}
\caption{}
\label{first}
\end{figure}
i.e. we take a projection $S^1\times S^1\to S^1$ and add a $1$-handle which only gets mapped to one half of $S^1$. This $f$ has 2 critical points, both of index 1.
With the loops in Figure \ref{first} we have $G=\pi_1(M)=\langle a_1,b_1,a_2,b_2\,|\,[a_1,b_1][a_2,b_2]=1\rangle$. The homomorphism $\xi=f_\#:\pi_1(M)\to\mathbb{Z}$ is then given by $\xi(a_1)=-1$ and all other generators are send to $0$.\\[0.2cm]
With the construction of \ref{CHtyTp} we get $N=S^1$ and $M_N$ as in Figure \ref{secon}. We can choose $N$ so that it is the image of the loop $b_1$, the basepoint being on the bottom.
Put an $f$-gradient $v$ on $M$ so that the trajectories of $-v$ starting and ending at critical points are as in Figure \ref{secon}.
\begin{figure}[ht]
\begin{center}
\epsfig{file=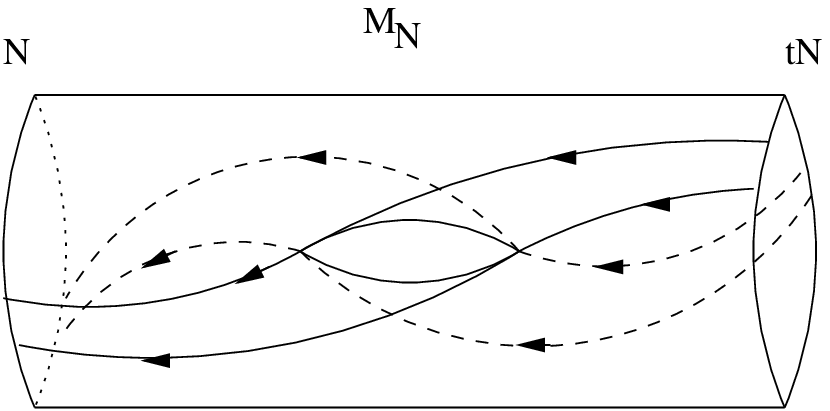}
\end{center}
\caption{}
\label{secon}
\end{figure}\\
We need $v$ to satisfy condition $(\mathfrak{C}')$. The Morse map $\psi:N\to\mathbb{R}$ is chosen as the height function, so we have a minimum and a maximum.
If the thickenings of the critical points on $N$ are chosen to fill about half of the circle it is clear that we can find a $v$ that satisfies $(\mathfrak{C}')$.
Now we can also get a $v\in\mathcal{G}_0(f)$ with trajectories as in Figure \ref{secon}.\\[0.2cm]
Let $\varphi(v):\widehat{\mathbb{Z}G}_\xi\otimes_{\mathbb{Z}G}C_\ast^\Delta(\tilde{M})\to C_\ast(f,v)$ be the chain equivalence from Theorem \ref{maint}. To calculate $\tau(\varphi(v))$ we have to look at the $1\times 1$ matrices $A_0$ and $A_1$ that come from the negative gradient descent.
All trajectories that start in $tN$ and are not drawn in Figure \ref{secon} flow to $N$ and cannot cross each other.
To calculate $A_0$ notice that trajectories starting in the lower half of $tN$ follow the loop that represents $a_2a_1$ up to conjugacy. Therefore $A_0=(a_2a_1)$.
The trajectories starting in the upper half of $tN$ and ending in the upper half of $N$ follow the loop $a_1$ up to conjugacy, so $A_1=(a_1)$.
Therefore 
\[\tau(\varphi(v))=\tau(1-a_1)-\tau(1-a_2a_1)=\tau((1-a_1)(1-a_2a_1)^{-1})\in\overline{K}_1^G(\widehat{\mathbb{Z}G}_\xi).\]
By Theorem \ref{maint} we get 
\[\eta(-v)=\varepsilon(\log((1-a_1)(1-a_2a_1)^{-1}))\in\widehat{\mathbb{Q}\Gamma}_\xi\]
and 
\[\zeta(-v)=(1-[a_1])(1-[a_2a_1])^{-1}\in\widehat{\mathbb{Q}H}_{\bar{\xi}}.\]
\label{examp}
\end{exam}
\begin{rema}\label{rm54}\em
Notice that in Figure \ref{secon} the unstable manifolds of $v$ intersect $tN$ in the upper half
of $tN$ while the stable manifolds intersect $N$ in the lower half of $N$. This allows $v$ to
satisfy $(\mathfrak{C}')$ with the Morse map $\psi$ on $N$. If we push the unstable manifolds
down and the stable manifolds up in Figure \ref{secon}, we get a different vector field $w$
which also satisfies $(\mathfrak{C}')$, but with the Morse map $-\psi$. So if we want to calculate
$\tau(\varphi(w))$ we have to interchange the roles of $A_0$ and $A_1$ which gives
\[\tau(\varphi(w))=\tau(1-a_2a_1)-\tau(1-a_1)=-\tau(\varphi(v))\in\overline{K}^G_1(\widehat{\mathbb{Z}G}
_\xi),\]
and
\[\zeta(-w)=(1-[a_1])^{-1}(1-[a_2a_1])\in\widehat{\mathbb{Q}H}_{\bar{\xi}}.\]
We can interpret this the following way: By looking at Figure \ref{secon} we can expect two closed
orbits, one on top of the cobordism, call it $\gamma_1$, and one on the bottom, call it $\gamma_2$.
The conjugacy class represented by $\gamma_1$ is the class of $a_1$ while $\gamma_2$ represents
the conjugacy class of $a_2a_1$. Now $\varepsilon(\gamma_1)=-1$ and $\varepsilon(\gamma_2)=1$
for the vector field $v$, but by passing to $w$ the unstable and stable manifolds move and the
signs change.
\end{rema}
\begin{appendix}
\section{The geometric chain homotopy equivalence}
The chain homotopy equivalence given by formula (\ref{formula}) has its counterpart in the exact
case. The purpose of this appendix is to describe the properties in that case. An alternative approach
can be found in Hutchings and Lee \cite[\S 2.3]{hutlee}, see also Schwarz \cite[\S 4.2]{schwar},
but since we need the torsion of the equivalence, we give full proofs.
\subsection{The relative Morse-Smale complex}
\label{app1}
Let $(W;M_0,M_1)$ be a compact cobordism, $f:W\to[a,b]$ a Morse function and $v$ an $f$-gradient
satisfying the transversality condition. A smooth triangulation $\Delta$ of $W$ is said to be \em
adjusted to \em $v$, if every $k$-simplex $\sigma$ intersects the unstable manifolds $D_R(p)$ transversely
for all critical points $p$ of index $\geq k$. In particular, if $p$ is a critical point of index
$k$, a $k$-simplex $\sigma$ intersects $D_R(p)$ in finitely many points. Using the orientations
we can assign to every such point a sign. Given a regular covering space $q:\tilde{W}\to W$ we can use
the covering transformation group $G$ and liftings of critical points and simplices to
assign an element $[\sigma:p]\in\mathbb{Z}G$ to the intersection and define a map \\
\begin{minipage}{3.5cm}\begin{eqnarray}\label{chheq}\\ \nonumber\end{eqnarray}
\end{minipage}\begin{minipage}{8cm}
\begin{eqnarray*}
\varphi:C^\Delta_\ast(\tilde{W},\tilde{M}_0)&\longrightarrow& C^{MS}_\ast(\tilde{W},\tilde{M}_0)\\
\sigma_k&\mapsto&\sum\limits_{p\in{\rm crit}_k(f)}[\sigma:p]\,p  \\
\end{eqnarray*}\end{minipage}

Here $C^{MS}_\ast(\tilde{W},\tilde{M}_0)$ is the Morse-Smale complex generated by the critical points
of $f$. For $A\subset W$ we denote $\tilde{A}=q^{-1}(A)$. Before we show the existence of adjusted
triangulations let us show that $\varphi$ is indeed a chain map.
\begin{lem}
\label{alem1}
$\varphi$ is a chain map.
\end{lem}
\begin{proof}
There exists a filtration $M_0=W_{-1}\subset W_0\subset\ldots\subset W_n=W$ of $W$ such that $W_i$
is a compact cobordism containing all critical points of index $\leq i$ and such that
$C^{MS}_k(\tilde{W},\tilde{M}_0)=H_k(\tilde{W}_k,\tilde{W}_{k-1})$ and the boundary homomorphism
comes from the long exact sequence of the triple $(\tilde{W}_k,\tilde{W}_{k-1},\tilde{W}_{k-2})$,
see Milnor \cite[\S 7]{milnhc}. Also, $C^\Delta_k(\tilde{W},\tilde{M}_0)=H_k(\tilde{W}^{(k)},
\tilde{W}^{(k-1)})$, where $W^{(k)}$ denotes the $k$-skeleton of the triangulation. A simplex
$\sigma_k\in C^\Delta_k(\tilde{W},\tilde{M}_0)$ is represented by a map $f_\sigma:(\Delta^k,
\partial\Delta^k)\to(\tilde{W}^{(k)},\tilde{W}^{(k-1)})$. Let $\Phi:\tilde{W}\times \mathbb{R}
\to\tilde{W}$ be induced by the flow of $-v$, where a flowline is supposed to stop once it hits
the boundary. For $t\in\mathbb{R}$ let $\Phi_t=\Phi(\cdot,t)$. Since $\Delta$ is adjusted to
$v$ there is a $t>0$ such that $\Phi_t\circ f_\sigma$ maps $\Delta^k$ to $\tilde{W}_k$ and
$\partial\Delta^k$ to $\tilde{W}_{k-1}$. It follows from intersection theory that
\[\varphi(\sigma)=(\Phi_t\circ f_\sigma)_\ast[\Delta^k]\in H_k(\tilde{W}_k,\tilde{W}_{k-1}).\]
Furthermore this does not depend on $t$ as long as $t$ is large enough. A diagram chase gives that
$\varphi$ is a chain map.
\end{proof}
Now we want to show the existence of adjusted triangulations. Let $\psi:W\to W$ be a diffeomorphism
homotopic to the identity and $\Delta$ a smooth triangulation of $W$. Then $\psi\Delta$ is the
triangulation of $W$ where simplices are composed with $\psi$. The corresponding chain complexes
can be identified by choosing a lifting $\tilde{\psi}:\tilde{W}\to\tilde{W}$.\\[0.2cm]
So let $\Delta$ be any smooth triangulation and $\psi_{-1}={\rm id}_W$. We can adjust $\psi_{-1}$
near the $0$-skeleton so that $0$-simplices intersect all unstable manifolds transversely. Since
the boundary of $W$ is transverse to the flow, we can leave it invariant. This way we get a
diffeomorphism $\psi_0$ isotopic to the identity. Now assume $\psi_{k-1}$ is isotopic to the
identity and every $j$-simplex of $\psi_{k-1}\Delta$ with $j\leq k-1$ intersects the unstable
manifolds transversely for critical points with index $\geq k-1$. We modify $\psi_{k-1}$ on the
$k$-skeleton so that $k$-simplices intersect $D_R(p)$ transversely for all $p$ with index $\geq k$.
Notice that for a $k$-simplex of $\psi_{k-1}\Delta$ this is already true for $\psi_{k-1}$ near the
boundary so we can leave the $(k-1)$-skeleton fixed. This way we obtain $\psi_k$ isotopic to the
identity and we can proceed by induction.\\[0.2cm]
Then $\psi_{n-1}\Delta$ is adjusted to $v$. Furthermore we can find an adjusted triangulation
$\psi\Delta$ with $\psi$ as close as we like to the identity. Moreover, compactness gives that
if $\Delta$ is adjusted to $v$, so is $\psi\Delta$ for every $\psi$ close enough to the identity.
On the other hand, given a triangulation $\Delta$ and two diffeomorphisms $\psi_1,\psi_2$ homotopic
to the identity such that $\psi_1\Delta$ and $\psi_2\Delta$ are adjusted to $v$, we get two chain
maps $\varphi_1$ and $\varphi_2$ which can be different.
\begin{lem}
\label{alem2}
The liftings can be chosen so that $\varphi_1$ and $\varphi_2$ are chain homotopic.
\end{lem}
\begin{proof} Let $H':W\times I\to W$ be a homotopy between $\varphi_1$ and $\varphi_2$. As above
we can change $H'$ to a homotopy $H:W\times I\to W$ between $\varphi_1$ and $\varphi_2$ such
that $H(\sigma\times I)$ intersects $D_R(p)$ transversely for all critical points $p$ with $
{\rm ind}p\geq k+1$, where $\sigma$ is a $k$-simplex. Describe the Morse-Smale complex as in the
proof of Lemma \ref{alem1}. Then we define $H_k:C^\Delta_k(\tilde{W},\tilde{M}_0)\to C^{MS}_{k+1}
(\tilde{W},\tilde{M}_0)$ by $H_k(\sigma)=(-1)^k(\Phi_t\circ\tilde{H})_\ast[\sigma\times I]\in
H_{k+1}(\tilde{W}_{k+1},\tilde{W}_k)$. Here $t>0$ is so large that $\Phi_t(\sigma\times\{0,1\}\cup
\partial\sigma\times I)\subset\tilde{W}_k$, $\Phi_t(\partial\sigma\times\{0,1\})\subset\tilde{W}_{k-1}$
and $\tilde{H}:\tilde{W}\times I\to \tilde{W}$ is a lifting of $H$. Use $\tilde{H}_0$ and $\tilde{H}_1$
to identify the triangulated chain complexes. Then $H_k$ is the desired
chain homotopy.
\end{proof}
Notice that for a $k$-simplex $\sigma$ and a disc $D_R(p)$ where ${\rm ind}\,p=k+1$ $H(\sigma\times I)
\cap D_R(p)$ is a finite set. So together with liftings and orientations we can write the chain
homotopy as
\[H_k(\sigma)=\sum_{p\in{\rm crit}_{k+1}(f)}[\sigma:p]\,p\hspace{1cm}\mbox{where }[\sigma:p]\in
\mathbb{Z}G,\]
which is independent of the filtration and only involves intersection numbers.
\begin{theo}
\label{sheqth}
Let $f:W\to [a,b]$ be a Morse function, $v$ an $f$-gradient satisfying the transversality condition,
$\Delta$ a triangulation adjusted to $v$ and $q:\tilde{W}\to W$ a regular covering space. Then $\varphi:C^{\Delta}_\ast(\tilde{W},\tilde{M}_0)
\to C_\ast^{MS}(\tilde{W},\tilde{M}_0)$ given by (\ref{chheq}) is a simple homotopy equivalence.
\end{theo}
\begin{proof}
Let $\Delta'$ be a subdivision of $\Delta$. If $\psi\Delta'$ is adjusted to $v$, so is $\psi\Delta$.
Moreover, the diagram
\[\begin{array}{rcl}C_\ast^{\psi\Delta}(\tilde{W},\tilde{M}_0)&\stackrel{\rm sd}{\longrightarrow}&
C_\ast^{\psi\Delta'}(\tilde{W},\tilde{M}_0)\\[0.2cm]
\varphi\searrow& &\swarrow\varphi\\[0.2cm]
\multicolumn{3}{c}{C_\ast^{MS}(\tilde{W},\tilde{M}_0)}\end{array}\]
commutes, where sd is subdivision, a simple homotopy equivalence. By Munkres \cite[\S 10]{munkre}
it is good enough to show the theorem for a special smooth triangulation.\\[0.2cm]
As in the proof of Lemma \ref{alem1} we have the filtration $M_0=W_{-1}\subset W_0\subset\ldots
\subset W_n=W$ where $W_i$ is a compact cobordism containing the critical points of index $\leq i$.
Choose a triangulation such that each $W_i$ is a subcomplex for all $-1\leq i \leq n$ and so that
for each critical point $p$ of index $i$ the disc $D_i(p)=D_L(p)\cap(W_i-{\rm int}\,W_{i-1})$ is a
subcomplex. We set for $0\leq k\leq n$ $C_\ast^{(k)}=C_\ast^{\Delta}(\tilde{W}_k,\tilde{M}_0)$.
The complex $D^{(k)}_\ast$ is given by
\[D_i^{(k)}=\left\{\begin{array}{cl}C^{MS}_i(\tilde{W},\tilde{M}_0)&i\leq k\\
0&\mbox{otherwise}\end{array}\right.\]
The chain map $\varphi$ induces maps $\varphi^{(k)}:C_\ast^{(k)}\to D_\ast^{(k)}$ and
$\varphi^{(k,k-1)}:C_\ast^{(k)}/C_\ast^{(k-1)}\to D_\ast^{(k)}/D_\ast^{(k-1)}$. Since the diagram
\[
\begin{array}{cclclclcc}
0& \longrightarrow & C_\ast^{(k-1)}& \longrightarrow& C_\ast^{(k)}&\longrightarrow& C_\ast^{(k)}/
C_\ast^{(k-1)}& \longrightarrow & 0\\[0.3cm]
& &\Big\downarrow \varphi^{(k-1)}& &\Big\downarrow\varphi^{(k)}& &\Big\downarrow\varphi^{(k,k-1)}\\[0.3cm]
0&\longrightarrow&D_\ast^{(k-1)}&\longrightarrow& D_\ast^{(k)}&\longrightarrow& D_\ast^{(k)}/
D_\ast^{(k-1)}&\longrightarrow& 0\end{array}\]
commutes, it suffices to show that each $\varphi^{(k,k-1)}$ is a simple homotopy equivalence to
finish the proof.\\[0.2cm]
Clearly $\varphi^{(k,k-1)}$ induces an isomorphism in homology, so it remains to show that it is
simple. We set $D_i=\bigcup\limits_{p\in{\rm crit}_i(f)}D_i(p)$. Then the inclusion $i:C_\ast
^\Delta(\tilde{W}_{i-1}\cup\tilde{D}_i,\tilde{W}_{i-1})\to C_\ast^\Delta(\tilde{W}_i,\tilde{W}_{i-1})$
is the inclusion of the core of the handles into the handles, hence a simple homotopy equivalence.
Now $\varphi^{(k,k-1)}\circ i$ is a simple homotopy equivalence by Cohen \cite[18.3]{cohen}, since
we can choose the lifts of $D_i$ so that the matrices representing $\varphi^{(k,k-1)}\circ i$ and
the boundary operators have only integer values. Therefore $\varphi^{(k,k-1)}$ is a simple homotopy
equivalence.
\end{proof}
\begin{rema}\em
Pajitnov \cite[Appendix A]{pajito} describes a simple homotopy equivalence \\$\psi:C_\ast^{MS}(\tilde{W},
\tilde{M}_0)\to C_\ast^\Delta(\tilde{W},\tilde{M}_0)$, where the triangulation is given by
\cite[Lm.A.8]{pajito}. The liftings can be chosen so that $\varphi\circ\psi$ is the identity on
the Morse-Smale complex, so $\varphi$ and $\psi$ are mutually inverse equivalences.
\end{rema}
\subsection{The absolute Morse-Smale complex}
\label{app2}
To calculate the absolute homology $H_\ast(\tilde{W})$ for a cobordism with nonempty boundary we
need more technicalities. So if $f:W\to[a,b]$ is a Morse function on a compact cobordism and $v$ an
$f$-gradient, Pajitnov \cite{pajirn,pajiov} defines sets $D_\delta({\rm ind}\leq i;v)$ and
$C_\delta({\rm ind}\leq i;v)$ for $\delta>0$ and $0\leq i\leq n$ which form filtrations of the
cobordism. To simplify the notation we denote them by $D_\delta^i(v)$ and $C_\delta^i(v)$. Pajitnov
\cite[Df.4.2]{pajiov} now defines a condition $(\mathfrak{C})$ on the vector field which also involves
Morse functions $\phi_i$ on $M_i$, $\phi_i$-gradients $u_i$ for $i=0,1$ and a $\delta>0$.\\[0.2cm]
Now we define a filtration of $W$ by $W_i=C^i_\delta(u_0)\cup D^i_\delta(v)$. It follows from the
methods of Pajitnov \cite[\S 5]{pajirn} that $C^{MS}_i(\tilde{W})=H_i(\tilde{W}_i,\tilde{W}_{i-1})$
gives a free $\mathbb{Z}G$ complex calculating $H_\ast(\tilde{W})$ with $C_i^{MS}(\tilde{W})=
C^{MS}_i(\tilde{M}_0)\oplus C_i^{MS}(\tilde{W},\tilde{M}_0)$. We want to explicitely describe
a chain homotopy equivalence between $C^\Delta_\ast(\tilde{W})$ and $C^{MS}_\ast(\tilde{W})$ based
on \ref{app1}.\\[0.2cm]
We need a triangulation $\Delta$ of $W$ with subtriangulations $\Delta_i$ of $M_i$ for $i=0,1$
with the properties 1.-4. described in Lemma \ref{trialem} below.
Notice that for $\varepsilon>\delta>0$ we have \\$C_\varepsilon^k(u_i)\subset C_\delta^k(u_i)$.\\[0.2cm]
So if $\sigma$ is a $k$-simplex of $\Delta$ we can flow it along $-v$ as in the proof of Lemma
\ref{alem1} into $W_k$. This induces a chain map $\varphi:C^\Delta_\ast(\tilde{W})\to C^{MS}_\ast
(\tilde{W})$. Also $C^\Delta_k(\tilde{W})=C^\Delta_k(\tilde{M}_0)\oplus C^\Delta_k(\tilde{W},\tilde{M}_0)$
and $C^{MS}_k(\tilde{W})=C^{MS}_k(\tilde{M}_0)\oplus C^{MS}_k(\tilde{W},\tilde{M}_0)$. In this
decomposition we have $\varphi=\left(\begin{array}{cc}\varphi_1&\ast\\0&\varphi_2\end{array}\right)$,
where $\varphi_1:C^\Delta_\ast(\tilde{M}_0)\to C^{MS}_\ast(\tilde{M}_0)$ and $\varphi_2:C^\Delta_\ast
(\tilde{W},\tilde{M}_0)\to C^{MS}_\ast(\tilde{W},\tilde{M}_0)$ are the simple homotopy equivalences
of \ref{app1}. Therefore $\varphi$ is also a simple homotopy equivalence.
\subsection{Existence of a nice triangulation}
\label{lastap}
Let $(W;M_0,M_1)$ be a compact cobordism, \\$f:W\to[-\frac{1}{2},n+\frac{1}{2}]$ an ordered Morse
function and $v$ an $f$-gradient satisfying the transversality condition and condition $(\mathfrak{C})$
of Pajitnov \cite[\S 4]{pajiov}. Let $\phi_i:M_i\to\mathbb{R}$, $u_i$ be the $\phi_i$-gradients and
$\delta>0$ given through condition $(\mathfrak{C})$.
\begin{lem}
\label{trialem}
There exists a triangulation $\Delta$ of W having $M_0\cup M_1$ as a subcomplex $\Delta_0\cup\Delta_1$
with the following properties:
\begin{enumerate}
\item $\Delta$ is adjusted to $v$.
\item For $i=0,1$ $\Delta_i$ is adjusted to $u_i$.
\item There is an $\varepsilon>\delta$ such that if $\sigma$ is a $k$-simplex in $\Delta_i$,
then $\sigma\subset C^k_\varepsilon(u_i)$.
\item There is an $\varepsilon>\delta$ such that if $x\in W^{(k)}$ and the trajectory of $-v$ starting
at $x$ ends in $fl(x)\in M_0$, then $fl(x)\in C_\varepsilon^k(u_0)$.
\end{enumerate}
\end{lem}
\begin{rema}\em Notice that all conditions in Lemma \ref{trialem} are open in the sense that if
$\Delta$ satisfies 1.-4., so does $\psi\Delta$, provided $\psi$ is close enough to the identity
in the smooth topology. Therefore a triangulation as needed in the proof of Proposition \ref{ptor}
also exists.
\end{rema}
\begin{proof}[Proof of Lemma \ref{trialem}]
For $i=0,1$ choose triangulations $\Delta_i$ of $M_i$ adjusted to $u_i$ and let \\$\Phi^i:M_i\times
\mathbb{R}\to M$ be the flow of $-u_i$. Then there is a $t>0$ such that 3.\ is satisfied for
$\Phi^i_t\Delta_i$. Extend $\Phi_t^0\Delta_0\cup\Phi^1_t\Delta_1$ to a triangulation $\Delta$ of $W$.
Choose a diffeomorphism $\psi$ so close to the identity that $\psi\Delta$ is adjusted to $v$ and
2.\ and 3.\ still hold. Modify this triangulation again so that if $\sigma$ is a $k$-simplex and
$x\in\sigma$ flows to $fl(x)\in M_0$, then $fl(x)\notin D_R(q;u_0)$ for critical points $q$ of
$\phi_0$ with ${\rm ind\,}q\geq k+1$. Notice that $D_R(q,u_0)$ is at most $n-k-2$ dimensional.
By making only very small changes we now have a triangulation
satisfying 1.-3.\ and by compactness the condition of 4.\ for some $\varepsilon>0$, but not
necessarily for $\varepsilon>\delta$. Rename this triangulation $\Delta$ and the two subcomplexes
of the boundary $\Delta_0$ and $\Delta_1$.\\[0.2cm]
Notice that condition 4.\ already holds for $\Delta_0$ and $\Delta_1$. This is trivial for $\Delta_0$
and for $\Delta_1$ it follows from condition $(\mathfrak{C})$. By continuity it also holds in a
small collar neighborhood of the boundary. We can think of this collar as $f^{-1}([-\frac{1}{2},
-\frac{1}{2}+\eta)\cup (n+\frac{1}{2}-\eta,n+\frac{1}{2}])$ for some $\frac{1}{2}>\eta>0$ and
assume that there are no 0-simplices in $f^{-1}((n+\frac{1}{2}-\eta,n+\frac{1}{2}))$. Let $\xi:W
\to[0,1]$ be a smooth function which is 1 outside of the collar and 0 in a smaller collar. Let
$\Phi:W\times\mathbb{R}\to W$ be the flow of the vector field $-\xi\cdot v$. There is a $T_1>0$
such that if $\sigma$ is a $k$-simplex in $\Delta$ which does not meet $M_1$, then $\Phi_t(\sigma)
\subset W_k=f^{-1}([-\frac{1}{2},k+\frac{1}{2}])$. Furthermore $\Phi_{T_1}\Delta$ satisfies
1.-3.\ and the same form of 4.\ as $\Delta$ does.\\[0.2cm]
Let $V_k=f^{-1}(\{k-\frac{1}{2}\})$ for $k=0,\ldots\hspace{-1pt},n$ and $U_k\subset f^{-1}[k-\frac{1}{2},
k))$ be diffeomorphic to $V_k\times[0,1]$ with $(x,0)$ corresponding to $x$ and $(x,t)$ lying on
the same trajectory of $v$. Let \\$X_k=W_k-({\rm int\,}U_k\cup W_{k-1})$, then $X_k$ is a compact
cobordism. By changing $f$ if necessary we can assume that $f|_{X_k}$ is a Morse function on this
cobordism.\\[0.2cm]
The Morse function $\phi_0:M_0\to\mathbb{R}$ is ordered so we can assume that $Y_k=\varphi_0^{-1}
((-\infty,k+\frac{1}{2}])$ gives a filtration with $D^k_\delta(u_0)\subset Y_k\subset C_\delta^k(u_0)$
for $k=0,\ldots\hspace{-1pt},n-1$. Let $\mu_k:M_0\to[0,1]$ be a smooth function with $\mu_k|_{Y_{k-1}}=0$
and $\mu_k|_{M_0-Y_k}=1$ and define $\nu_k=\mu_k\cdot u_0$ and let $\Lambda_k:M_0\times\mathbb{R}
\to M_0$ be the flow of $-\nu_k$.\\[0.2cm]
Define $\Theta_k:X_k\times\mathbb{R}\to X_k$ as follows. Let $(x,t)\in X_k\times\mathbb{R}$.
If $x\in\bigcup\limits_{{\rm ind}p\leq k}D_R(p;v)$, let $\Theta_k(x,t)=x$. If not, the trajectory
of $-v$ reaches $fl(x)\in M_0$. Let $fl_t(x)=\Lambda_k(fl(x),t)$. Now $fl_t(x)$ flows along $v$
all the way back to $y\in f^{-1}(\{f(x)\})$. For if not, we have $fl_t(x)\in\bigcup\limits_{{\rm ind}
p\leq k}D_L(p;v)$. By $(\mathfrak{C})$ we have $fl_t(x)\in Y_{k-1}$. But then $fl(x)=fl_t(x)$
because $\nu_k(fl_t(x))=0$ and then $fl_t(x)$ flows all the way back to $x$. This and the implicit
function theorem give that $\Theta_k(x,t)=y$ is smooth on $X_k\times\mathbb{R}$. Similarly,
a smooth homotopy $\lambda_k$ between $\mu_k$ and $\mu_{k-1}$ with $\lambda_k|_{(M_0-Y_k\times[0,1]}
=1$ and $\lambda_k|_{Y_{k-2}\times[0,1]}=0$ defines a smooth map $\Theta_k':U_k\times\mathbb{R}\to U_k$
such that all $\Theta_k$ and $\Theta_k'$ define an isotopy $\Theta:W\times\mathbb{R}\to W$ with the
following property : if $x\in W_k-\bigcup\limits_{{\rm ind}p\leq k}D_R(p;v)$ with $fl(x)\in M_0-\bigcup
\limits_{{\rm ind}q\geq k+1} D_R(q,u_0)$, then there is a $t_x>0$ such that $\Theta(x,s)\in W_k
-\bigcup\limits_{{\rm ind}p\leq k}D_R(p)$ and $fl(\Theta(x,s))\in C_\delta^k(u_0)$ for all $s\geq t_x$.
By compactness there is a $T_2>0$ such that $\Theta_{T_2}\Phi_{T_1}\Delta$ satisfies the Lemma.
\end{proof}
\end{appendix}

\end{document}